\documentclass[a4paper, 11pt]{amsart}
\usepackage[utf8]{inputenc}
\usepackage{graphicx}
\usepackage{stmaryrd}
\usepackage{xcolor}
\usepackage{verbatim}

\usepackage{bbold}
\usepackage[margin=3cm]{geometry}
\usepackage{float}
\usepackage[utf8]{inputenc}
\usepackage{amsthm}
\usepackage{mathtools}
\usepackage{amssymb}
\usepackage{amsmath}
\usepackage{blkarray}
\usepackage{enumitem}

\newcommand{\N}{\mathbb N}

\newcommand{\C}{\mathbb C}

\newcommand{\Hom}{\normalfont\text{Hom}}
\newcommand{\ev}{\normalfont\text{ev}}

\newcommand{\op}{\normalfont\text{op}}
\newcommand{\defn}{\text{def}}

\newcommand{\V}{\mathcal{V}}
\newcommand{\Lagr}{\mathcal{L}}
\newcommand{\E}{\mathcal{E}}
\newcommand{\R}{\mathcal{R}}
\newcommand{\A}{\mathbb{A}}

\newcommand{\B}{\mathbb{B}}

\newcommand{\sS}{\mathbf{sSet}}
\newcommand{\sE}{\mathbf{s\E}}
\newcommand{\s}{\mathbf{Set}}

\newcommand{\Gpd}{\mathbf{Gpd}}
\newcommand{\GpdE}{\mathbf{Gpd(\E)}}
\newcommand{\Cat}{\mathbf{Cat}}
\newcommand{\CatE}{\mathbf{Cat(\E)}}
\newcommand{\X}{\mathbb{X}}
\newcommand{\Y}{\mathbb{Y}}
\newcommand{\Par}{\mathcal{P}}
\newcommand{\K}{\mathcal{K}}

\newcommand{\lifts}{\boxslash}
\newcommand{\elifts}{\underline{\boxslash}}

\newcommand{\splt}{\text{split}}

\newcommand{\Weq}{\mathbf{W}}
\newcommand{\Cof}{\mathbf{C}}

\newcommand{\gencof}{\mathbf{I}}
\newcommand{\gentrivcof}{\mathbf{J}}
\newcommand{\1}{\mathbf{1}}
\newcommand{\TE}{\underline{\1}}
\newcommand{\TF}{\mathbf{TF}}
\newcommand{\F}{\mathbf{F}}
\newcommand{\TC}{\mathbf{TC}}
\newcommand{\Map}{\mathbf{Map}}
\newcommand{\iso}{\mathbf{Iso}}

\newcommand{\disc}{\text{disc}}

\usepackage{color}
\usepackage{hyperref}

\usepackage[    noabbrev,
        capitalise]{cleveref}

\theoremstyle{plain}
\newtheorem{thm}{Theorem}[section]
\newtheorem{lem}[thm]{Lemma}
\newtheorem{cor}[thm]{Corollary}

\newtheorem{prop}[thm]{Proposition}

\theoremstyle{definition}

\newtheorem{define}[thm]{Definition}

\newtheorem{remark}[thm]{Remark}

\crefname{prop}{Proposition}{Propositions}
\crefname{thm}{Theorem}{Theorems}
\crefname{lem}{Lemma}{Lemmas}

\usepackage{tikz-cd}
\usepackage{tikz}
\usetikzlibrary{arrows}
\usepackage{pgfplots}
\pgfplotsset{compat=1.18}
\usetikzlibrary{positioning, calc}
\usetikzlibrary{decorations.markings,arrows.meta,bending}
\tikzset{-->--/.style={decoration={markings, 
			mark=at position #1 with {\arrow[line width=2pt]{>}}},postaction={decorate}}}

\NewDocumentCommand\Cycle{O{} m m m O{} m}{%
	\draw[#1](#2.{#3+asin(#6/(#4*1.41))}) arc (180+#3-45:180+#3-45-270:#6/2) #5;
}

\tikzset{
	partial ellipse/.style args={#1:#2:#3}{
		insert path={+ (#1:#3) arc (#1:#2:#3)}
	}
}
\tikzset{->-/.style={decoration={
			markings,
			mark=at position .5 with {\arrow{>}}},postaction={decorate}}}

\usepackage{tkz-euclide}
\usepackage{tikz-3dplot}
\tdplotsetmaincoords{60}{60}

\usepackage{stackengine, graphicx}

\begin{document}

\title{The algebraic internal groupoid model of Martin-L\"{o}f type theory} 


\author{Calum Hughes} 

\begin{abstract}
  We extend the model structure on the category $\CatE$ of internal categories studied by Everaert, Kieboom and Van der Linden to an algebraic model structure. Moreover, we show that it restricts to the category of internal groupoids. We show that in this case, the algebraic weak factorisation system that consists of the algebraic trivial cofibrations and algebraic fibrations forms a model of Martin-L\"{o}f type theory. Taking $\E = \s$ and forgetting the algebraic structure, this recovers Hofmann and Streicher's groupoid model of Martin-L\"{o}f type theory. Finally, we are able to provide axioms on a $(2,1)$-category which ensure that it gives an algebraic model of Martin-L\"{o}f type theory. To do this, we give necessary and sufficient axioms on a $2$-category $\K$ such that $\K\simeq \CatE$ in which $\E$ is a locally cartesian closed locos with coequalisers, a result which we believe is of independent interest. 
\end{abstract}

\maketitle
\tableofcontents

\section{Introduction}

\subsection{Context and motivation}

The groupoid model of Martin-L\"{o}f type theory (MLTT) was introduced by Hofmann and Streicher in 1998 \cite{hofmann1998groupoid}, proving that the uniqueness of identity proofs property cannot be derived in this type theory and thus revealing a deep connection between logic, (higher) category theory and homotopy theory (see also the work of Lamarche \cite{Lamarche1991Proposal}). This key observation led to the development of homotopy type theory and univalent foundations \cite{program2013homotopy}. Since then, many different categorical models of (homotopy) type theory have been studied with a particular focus on those which have an abstract homotopy theory \cite{voevodsky2006very,  bezem2014model, riehl2017type, berg2018path, berg2020univalent,awodey2023cartesian}. 

One of the earliest and most effective ways to abstractly provide a category with an abstract homotopy theory is given by Quillen's notion of a model structure \cite{quillen1967homotopical}. Loosely, this is a category with three distinguished classes of maps called weak equivalences, fibrations and cofibrations. From a categorical point of view these suffer from certain defects, for example cofibrations (resp. fibrations) are not closed under colimits (resp. limits) in the arrow category. Riehl's notion of an algebraic model category \cite{riehl2011algebraic} builds on the notion of an algebraic weak factorisation system \cite{grandis2006natural, garner2009understanding} and provides a solution to this by considering cofibrations and fibrations as retracts of coalgebras for a comonad and algebras for a monad respectively, since (co)algebras for a (co)monad are closed under (co)limits.  The philosophy is that by keeping track of algebraic data, certain calculations become canonical and therefore more categorically well-behaved. 

An example of key interest to us is given by the classical model structure on $\Cat$, the category of small categories and functors and the restriction of this to the category of small groupoids, $\Gpd$ \cite{anderson1978fibrations, rezk1996model}. Here, the weak equivalences are equivalences of categories; fibrations are isofibrations; cofibrations are functors which are injective on objects. Because this model structures is cofibrantly generated, it underlies an algebraic model structure by a process described in \cite{riehl2011algebraic}, although the explicit details of this for the classical model structure on $\Cat$ has not been worked out until now. 

Hofmann and Streicher's groupoid model of MLTT lives in the context of the homotopy theory given by the classical model structure on $\Gpd$; dependent types are modelled by isofibrations between groupoids. In fact, there are a couple of coherence issues with this model; substitution and identity types are defined only up to isomorphism rather than equality. One suggested approach to fixing this problem is using precisely the algebraic perspective described above. By keeping track of algebraic data there are canonical choices for substitutions and identity types in the model. To this end, Gambino and Larrea introduce the notion of a type theoretic algebraic weak factorisation system \cite{gambino2023models}, which is an algebraic weak factorisation system with extra structure.

In this paper, we shall extend these ideas to the setting of internal categories and internal groupoids. These notions have recently witnessed an increase in interest partially due to their connection with double category theory, $2$-dimensional universal algebra and alternative foundations of mathematics \cite{johnstone2002sketches, everaert2005model, joyal2006strong, fiore2008model, bourke2010codescent, niefield2019internal, Hughes2024Colimits, Hughes2024Elementary}. The idea is to replace the \emph{sets} of objects and morphisms of a small category with \emph{objects} of a category $\E$ and seeing how much category theory can be internalised to this setting. Categories internal to $\E$ form a category which we denote $\CatE,$ and we denote its subcategory of internal groupoids by $\GpdE$. The homotopy theory of $\CatE$ has been studied in \cite{everaert2005model}. In \cite[\S 7]{everaert2005model}, a model structure is considered in which the fibrations are the representable isofibrations and the weak equivalences are the representable equivalences of categories. We call this the \emph{natural model structure} on $\CatE$. Note that this is not in general the same as the model structure on internal categories considered in \cite{joyal2006strong}. The goal of this paper is to extend this to an algebraic model structure and use it to obtain models of MLTT. 

\subsection{Main results}

This paper makes two main contributions. The first contribution is to prove some nice properties of the natural model structure on $\CatE$. We show that it is a cofibrantly generated, monoidal model structure which underlies an algebraic model structure (\cref{thm: internal analogue}) and give explicit descriptions of the algebraic structure. As an application, we conclude that for an internal monoidal category $\mathbb{M}$, there is a cofibrantly generated monoidal model structure on the category of $\mathbb{M}$-modules. This is \cref{cor: model structure on modules}.

The second contribution of this paper is that by restricting to $\GpdE$, we obtain an algebraic model of MLTT. We show that the algebraic weak factorisation system consisting of maps with a trivial cofibration structure and maps with a fibration structure forms a type theoretic algebraic weak factorisation system in the sense of \cite{gambino2023models}. This is \cref{thm: ttawfs}. 

The main benefit of this approach is that we can apply these results to a variety of different examples. In some cases, this gives an improved perspective on results that have been previously considered in the literature; in other cases, it gives novel results. For example, we are able to apply our result to groupoids internal to $\s$, the category of (modest) assemblies and the effective topos, and $[\C, \s]$ giving us constructive versions of classical results, realisability models of type theory, and an indexed version of the groupoid model of type theory respectively.

\subsection{Related work}

Algebraic groupoid models of MLTT are also considered in \cite[Theorem 5.5]{gambino2023models}. Our paper works in the more general setting groupoids internal to some category $\E$, and places this in the larger structure of an algebraic model structure. Taking $\E=\s$, the precise relationship between these models is explained in \cref{subsubsec: set}. 

We note that Anthony Agwu is studying models of type theories in groupoids internal to the effective topos and a model structure using different methodology. We believe that the model being considered is for \emph{normal} isofibrations in contrast to ours, which deals with \emph{cloven} isofibrations. As a result, both the model structure and the model of type theory will be different. The method of looking at groupoids internal to the effective topos as a model of type theory is also related to the work of both Awodey and Emmenegger \cite{awodey2025effective} and Speight, who are both considering $2$-dimensional forms of the effective topos.

\subsection{Outline}

\cref{sec: prelim} establishes the tools which we will exploit in our proofs including an enrichment of $\CatE$ over $(\E, \times, \mathbf{1})$ which allows us to run the enriched small object argument on a pair of classes of maps in $\E$. In \cref{sec: cofib gen}, we show that these form the generated (trivial) cofibrations of the model structure of \cite[\S 7]{everaert2005model}, recalled in \cref{subsec: nat mod struc}; whence the model structure is cofibrantly generated. Using this, we can easily show that the model structure is symmetric monoidal and admits an algebraic weak model structure. This is recorded in \cref{thm: internal analogue}. We show that this restricts to $\GpdE$ in \cref{cor: algebraic model structure on GpdE}.  

In \cref{sec: Algebraic aspects}, we work through the details of the algebraic structure explicitly which allows us to show in \cref{sec: TTAWFS} that for $\E$ a locally cartesian closed locos with coequalisers, the data of the algebraic model structure on $\GpdE$ has the extra structure of a type theoretic algebraic weak factorisation system, in the sense of Definition 4.10 of \cite{gambino2023models}. This is \cref{thm: ttawfs}, and proves that we have an internal groupoid model of MLTT for any locally cartesian closed locos with coequalisers.

We conclude with \cref{sec: examples}, which applies the results of the paper to various examples giving both novel results and extensions of previously known results. 

\subsection{Acknowledgements}

The support of the Dame Kathleen Ollerenshaw PhD studentship is gratefully acknowledged. I would also like to thank Nicola Gambino and Adrian Miranda for their advice while preparing this paper, and to Sam Speight for helpful comments on the first draft.

 \section{Preliminaries} \label{sec: prelim}

Throughout this paper, we will assume that $\E$ is a lextensive cartesian closed category with pullback stable coequalisers in which the forgetful functor $\mathcal{U}: \CatE \to \mathbf{Gph}(\E)$ has left adjoint. These assumptions ensure that $\CatE$ is finitely cocomplete \cite{Hughes2024Colimits}. 

Examples of such a setting are given by any any locally cartesian closed locos with coequalisers such as an arithmetic $\Pi$-pretopos and an elementary topos with a natural numbers object. These are key motivational examples for us. Note that we do not require regularity, exactness nor a subobject classifier. \cref{sec: examples} gives more detailed descriptions of suitable $\E$.

\begin{remark}
\label{rem: locallyfp}
    We could instead have chosen $\E$ to be a lextensive category that is locally finitely presentable and cartesian closed, such as $\Cat$, $\mathbf{Gph}$ and $\mathbf{Ab}$--- this would also ensure that $\CatE$ is finitely cocomplete \cite{Hughes2024Colimits}. Hence, the existence of an algebraic model structure can be applied instead when $\E$ is a lextensive category that is locally finitely presentable and cartesian closed. However, being locally finitely presentable is rather restrictive and furthermore not an elementary condition. Moreover, when modelling type theory, we want isofibrations of internal groupoids to be exponentiable; in fact, we show in \cref{lem: lcc} that isofibrations in $\GpdE$ are exponentiable if and only if $\E$ is locally cartesian closed. As such, for type theoretic aspects of this work, we cannot take $\E$ to be either $\Cat$ or $\mathbf{Ab}$ due to their lack of local cartesian closure.   
\end{remark}

 \subsection{Notation}

We adopt the notation for internal category theory established in \cite[\S 2]{Hughes2024Elementary} with one notable exception: since everything we deal with in this paper is $1$-categorical, we denote the $1$-category of internal categories and internal functors by $\CatE$ as opposed to $\CatE_1$, which is used to distinguish between this and the $2$-category $\CatE$ in that paper. 

We denote by $\Delta_{\le 3}$ the category of non-empty ordinals up to and including the ordinal with four elements. 

The following categories will be useful.
\[\mathbf{1}  := \bullet \]
	\[
	\mathbf{2}  := \begin{tikzcd}
		\bullet \arrow[r] & \bullet 
	\end{tikzcd}\]
\[\mathcal{P}  := \begin{tikzcd}
	\bullet \arrow[r, bend left] \arrow[r, bend right] & \bullet 
\end{tikzcd} 
\]
\[\mathcal{I} : = \begin{tikzcd}
    \bullet \arrow[r, bend left] & \bullet \arrow[l, bend left]
\end{tikzcd}\]

For a $\mathcal{V}$-enriched category $\mathcal{C}$ and $X,Y \in \text{Ob}(\V)$, we denote its hom-object by $\Hom_{\V}(X,Y)$. We note that $\CatE$ is $\Cat$-enriched and make use of this in defining certain concepts representably with respect to hom-categories. 

For $E \in \E$ and $\X \in \CatE$, we write $E \times \CatE : = \disc(E) \times \CatE$ and $\X^E : = \X ^{\disc(E)}.$

For $(T, \eta)$ a pointed endofunctor (resp. $\mathbb{T}$ a monad), we denote its category of algebras by $(T, \eta)\text{-}\mathbf{Alg}$ (resp. $\mathbb{T}\text{-}\mathbf{Alg}).$ 

For $f,g \in \Cat^{\mathbf{2}}$, we distinguish between $(u,v): f \to g$ and $\alpha: f \Rightarrow g$; the first is a commutative square from $f$ to $g$ i.e. $u: \text{dom}(f) \to \text{dom}(g)$ and $v: \text{cod}(f) \to \text{cod}(g)$ such that $gu = vf$. This is a morphism in the category $\CatE^{\mathbf{2}}$. The latter is an internal natural transformation.

 \subsection{Enrichment over $\E$}\label{subsec: enriched}
	
We define an enrichment of $\CatE$ over $(\E, \times, \mathbf{1})$. The functor $\Hom_{\s}(\TE, -): \E \to \s$ has a partial left adjoint, defined by by mapping any finite set $S$ to \begin{equation*}
		\label{eq:underline}
		\underline{S} := S \cdot \TE = \coprod_{s \in S} \TE. \end{equation*}
	This can also be described as a (genuine) functor $\underline{(-)}: \mathbf{FinSet} \to \E$. 
	
    We extend this levelwise to become a partial functor $\underline{(-)
}: \Cat \to \CatE$, defined on \emph{finite categories}, equivalently a genuine functor $\underline{(-)}: \Cat(\mathbf{FinSet}) \to \CatE$. This is well-defined by extensivity of $\E$, so that in $\E$ the coproduct commutes with finite limits \cite{carboni1993introduction}; in particular, for $\C \in \Cat$, we have a well-defined composition operation
	
	\[\underline{m}: \underline{\C_1} \times_{\underline{\C_0}}\underline{\C_1} \to \underline{\C_1}.\]

 Hence we have a partial adjunction defined on the class of small categories with finite set of objects and finite set of morphisms. \begin{equation*}
			\begin{tikzcd}[ampersand replacement=\&, column sep=small]
				\Hom_{\Cat}(\TE, -) \vdash \underline{(-)}: \CatE \arrow[r, shift right =1] \arrow[r, leftharpoonup, shift left = 1] \& \Cat 
			\end{tikzcd}
\end{equation*}

\begin{remark}\label{rem: underline sets}
    For a set of finite categories $\mathbf{X} : = \{\C_0, \C_1,..., \C_n\}$ we define the set of internal categories $\underline{\mathbf{X}} : = \{ \underline{\C_0}, \underline{\C_1},..., \underline{\C_n}\}.$
\end{remark}

It is often useful to consider internal categories as simplicial objects. In our setting, since we do not want to assume countable limits, we consider internal categories as simplicial objects truncated at the fourth level. We can do this without loss of any information about the internal category because the highest data that internal categories have is associativity, which occurs at the fourth level of the simplicial nerve.  

\begin{lem}
    The truncated internal nerve functor $N: \CatE \hookrightarrow \left[\Delta_{\le 3}^{\op}, \E\right]$ is fully faithful.
\end{lem}

\begin{proof}
This is immediate from the definition of internal functor and the definition of morphism in $\left[\Delta_{\le 3}^{\op}, \E\right]$.
\end{proof}

Consider $\Delta_{\le 3}^{\op}$ with the free $\E$-enrichment as described in Example 3.4.4 of \cite{riehl2014categorical} in which $\Hom_{\E}\left([n],[m]\right) := \underline{\Hom_{\s}\left([n],[m]\right)}$ which is well-defined since $\Hom_{\s}\left([n], [m]\right)$ is finite. 

Let $\X, \Y$ be internal categories. Then $N\X$ and $N\Y$ are $\E$-valued presheaves. We define the following:
	
	\[\Hom_{\E}(\X,\Y) : = \int_{[n] \in \Delta_{\le 3}^{\op}} NY_n^{NX_n}. \]

 This is the hom-object of $\E$-natural transformations in $[\Delta_{\le 3}^{\op}, \E]$ as described in \cite[Digression 7.3.1]{riehl2014categorical}. By fully-faithfullness of the truncated internal nerve, this provides an enrichment on $\CatE.$

 \begin{remark}
     \label{rem: alternative description of enrichment}
     We give an equivalent description of the above enrichment. Recall from Proposition 3.2.1 of \cite{Hughes2024Elementary} that if $\E$ has finite limits and is cartesian closed, then $\CatE$ is cartesian closed with internal hom calculated in $\sE$, i.e. $\CatE$ is an exponential ideal of $\sE$. This factors through the truncated simplicial objects, $\left[\Delta_{\le 3}^{\op}, \E\right]$. Explicitly, for $\X, \Y \in \CatE,$ we have

     \begin{equation*}
         \underline{\mathbf{[}\X, \Y \mathbf{]}}  : = [N\X, N\Y] \in \CatE
     \end{equation*}

     which is defined to be the functor $\Delta_{\le 3}^{op} \to \E$ given by

     $$k \mapsto \displaystyle{\int\limits_{[n] \in \Delta_{\le 3}^{\op}}{Y_n}^{X_n \times \Delta_{\le 3}^{op} \left([n],[k]\right)}}.$$

     In particular, we have

    \begin{align*}
         \underline{\mathbf{[}\X, \Y \mathbf{]}}_0 & : = [N\X, N\Y]([0]) \\
         & = \displaystyle{\int\limits_{[n] \in \Delta_{\le 3}^{\op}} {NY_n}^{NX_n \times \Delta_{\le 3}^{op} \left([n],[0]\right)}} \\
         & = \displaystyle{\int\limits_{[n] \in \Delta_{\le 3}^{\op}}  {NY_n}^{NX_n \times \Delta_{\le3} \left([0],[n]\right)}} \\
         & \cong \displaystyle{\int\limits_{[n] \in \Delta_{\le 3}^{\op}}  {NY_n}^{NX_n \times \TE}} \\
         & \cong \displaystyle{\int\limits_{[n] \in \Delta_{\le 3}^{\op}}{NY_n}^{NX_n}} 
    \end{align*}

    So we equivalently have $\Hom_{\E}(\X, \Y) =  \underline{\mathbf{[}\X, \Y \mathbf{]}}_0.$ This different perspective will also be useful. 
 \end{remark}
	
The following follows from a strong form of the enriched Yoneda lemma \cite[Theorem 2.4]{kelly1982basic}.

\begin{prop}\label{prop: enriched yoneda}
     Let $\E$ be a cartesian closed category with finite limits. Let $X: \Delta_{\le 3} \to \E$. For all $[k] \in \Delta_{\le 3}$ we have that:

    \begin{equation*}
       X[k] \cong \int_{ [n] \in \Delta_{\le 3}} X_n^{\Delta([k],[n])}. 
    \end{equation*}
\end{prop}

	\begin{lem}
		\label{lem:yoneda}
		Let $\X$ be an internal category. We can calculate the following:
		\begin{enumerate}
			\item  $\Hom_{\E}(\underline{\emptyset}, \X) \cong \TE$.
			\item $\Hom_{\E}(\TE, \X) \cong X_0$.
			\item $\Hom_{\E}(\underline{\mathbf{2}}, \X) \cong X_1$.
			\item  $\Hom_{\E}(\TE + \TE, \X) \cong X_0\times X_0$.
   \item $\Hom_{\E}(\X, \TE) \cong \TE.$
			
		\end{enumerate}
	\end{lem}
	
	\begin{proof}
		(1) follows from the observation that $\underline{\emptyset}$ is the initial object in $\E$, and so for any $E \in \E$ we have $E^{\underline{\emptyset}}=\TE$. For (2) We have that $\underline{\Delta}_{\le3}^{\op} ([0],-) \cong N\TE$. Therefore, by \cref{prop: enriched yoneda}, we have $$\Hom_{\E}(\TE,\X) =_{\defn} [N\TE,N\X]_0 \cong [\underline{\Delta}_{\le3}([0],-),N\X] \cong N\X([0]) = X_0.$$ The proof for (3) is similar. For (4), we note that for enriched hom functors, colimits in the first variable can be taken out as limits \cite[\S 3.2]{kelly1982basic} and the result follows. For (5), note that $\TE^E \cong \TE$ for any $E \in \E.$
	\end{proof}
	
	From the internal hom, we can define an evaluation.
	\begin{define}
		We define a partial 2-variable functor $\ev_{(-)}(-): \Cat^{\op} \times \CatE \rightharpoonup \E.$ For $\mathcal{A}$ a finite category and $\X \in \CatE$:
		\[\ev_{\mathcal{A}}(\X): = \Hom_{\E}(\underline{\mathcal{A}}, \X).\]
		Fixing $\X$, we get a partial functor $\X(-) : = \ev_{(-)}(\X) : \Cat^{\op} \rightharpoonup \E.$

	\end{define}
	
	\begin{remark}
		\label{rem:eval}
		Note that \cref{lem:yoneda} then implies that for $\X \in \CatE$, we have:
		
		\begin{align*}
			\X(\emptyset) & \cong \TE \\
			\X(\mathbf{1}) & \cong X_0\\
			\X(\mathbf{2}) & \cong X_1\\
			\X(\mathbf{1}+\mathbf{1}) & \cong X_0 \times  X_0.
		\end{align*}
		
		This is an internal counterpart to statements in the ordinary setting such as maps from $\mathbf{1} $ to $\X$ are in bijection with its objects, $X_0$, which can be proven in $\Cat$ using the Yoneda lemma.

	\end{remark}

 Recall the functor $\iso: \CatE \to \GpdE$, the right adjoint to the inclusion $\GpdE \to \CatE$ (cf. \cite{bunge1979stacks}). We also have the following:

 \begin{lem}
 \label{lem: ev at I is iso}
     Let $\X \in \CatE_0$. We have $\X(\mathcal{I}) \cong \iso(\X)_1.$
 \end{lem}

 \begin{proof}
        Firstly, \cref{rem: alternative description of enrichment} tells us that $\X(\mathcal{I}) = \left(\X^{\underline{\mathcal{I}}}\right)_0$. The discussion of the construction of the power of $\X$ by the ordinary category $\mathcal{I}$ is given in section 3 of \cite{everaert2005model}, and has $\left(\X^{\mathcal{I}}\right)_0 :=\iso(\X)_1.$ Note that this construction is a groupoidal version of the power of an internal category by the category $\mathbf{2}$ discussed in section 3.1 of \cite{Hughes2024Elementary} and by the same argument given there, it follows that $\X^{\underline{\mathcal{I}}}$ has the universal property of $\X^{\mathcal{I}}$. We conclude that $\X(\mathcal{I}) \cong \iso(\X)_1.$
 \end{proof}

 We can also describe the object of pairs of parallel arrows in an internal category. Note that in $\Cat$, the category $\Par$ can constructed by gluing two walking arrows together, as in the pushout below:

    \begin{equation}\label{eq: pushout for P}
        \begin{tikzcd}
            \mathbf{1} + \mathbf{1} \arrow[d] \arrow[r] & \mathbf{2} \arrow[d] \\
            \mathbf{2} \arrow[r] & \Par \arrow[ul, phantom, "\ulcorner", very near start]
        \end{tikzcd}
    \end{equation}

 \begin{lem}
 \label{lem: parallel}
     We have a pullback square

     \begin{equation*}
         \begin{tikzcd}
             \X(\Par) \arrow[dr, phantom, "\lrcorner", very near start] \arrow[r] \arrow[d] & X_1 \arrow[d, "(d_1{,}d_0)"] \\
             X_1 \arrow[r, "(d_1{,}d_0)"'] & X_0 \times X_0
         \end{tikzcd}
     \end{equation*}
 \end{lem}

 \begin{proof}
     First, we apply $\underline{(-)}$ to the defining pushout square of $\Par$; this functor is (partial) left adjoint so preserves colimits of finite categories and results in a pushout square in $\E$. Now, some simple calculus using cartesian closedness and lextensivity of $\E$ shows that $\Hom_{\E}(-, \X)$ turns pushouts into pullbacks. Hence, we get the above pullback square. 
 \end{proof}

	This enrichment of $\CatE$ in $(\E, \times, \TE)$ makes it copowered and powered by $\E.$ Recall that this means for all $E \in \E$ and $\X \in \CatE,$ there exists objects $E \odot \X$, (the \emph{copower} of $\X$ by $E$) and $E \pitchfork \X$ (the \emph{power} of $\X$ by $E$) such that for any $\Y \in \CatE$ we have
		\[\Hom_{\E}(E \odot \X, \Y) \cong \Hom_{\E}(E, \Hom_{\E}(\X,\Y)) \cong \Hom_{\E}\left(\X, E \pitchfork \Y\right).\]

        Note that this condition is also called being tensored and cotensored over $\E$.

	\begin{prop}
 \label{prop: CatE tensored and cotensored}
		For any $E \in \E$ and $\X \in \CatE$, the copower of $\X$ by $E$ is given by $E \times \X$. Dually, the power of $\X$ by $\E$ is given by $\X^E$.
		
	\end{prop}
	
	\begin{proof}
		\begin{align*}
			\Hom_{\E}(E, \Hom_{\E}(\X,\Y))& := \Hom_{\E}\left(E, \int_{[n] \in \Delta_{\le3}^{\op}} NY_n^{NX_n} \right) & \\
			& \cong \int_{[n] \in \Delta_{\le3}^{\op}}\Hom_{\E}\left(E,  NY_n^{NX_n}\right) & \text{ as $\Hom_{\E}(E,-)$ is $\E$-continuous,}\\
			& \cong \int_{[n] \in \Delta_{\le3}^{\op}}\Hom_{\E}\left(E\times NX_n,  NY_n\right) & \text{tensor-hom adjunction,} \\
			& \cong \int_{[n] \in \Delta_{\le3}^{\op}} NY_n ^{E\times NX_n} & \text{ definition of enrichment of $\E$ over itself} \\
			& \cong \int_{[n] \in \Delta_{\le3}^{\op}} NY_n ^{(NE\times NX)_n} & \text{ definition of tensoring and nerve,} \\
			& \cong \int_{[n] \in \Delta_{\le3}^{\op}} NY_n ^{N(E\times X)_n} &  \\
			& =: \Hom_{\E}(E \times \X, \Y)  & \text{ as required.}\\
		\end{align*}
		Similarly:
		
		\begin{align*}
			\Hom_{\E}(E, \Hom_{\E}(\X,\Y))
			& \cong \int_{[n] \in \Delta_{\le3}^{\op}} NY_n ^{NE_n\times NX_n} & \\
			&  \cong \int_{[n] \in \Delta_{\le3}^{\op}} \left(NY_n ^{NE_n}\right)^{NX_n} & \\
			&\cong \int_{[n] \in \Delta_{\le3}^{\op}} N\left(Y ^{E}\right)_n^{NX_n} & \\
			& =: \Hom_{\E}\left(\X, \Y^E\right)  & \text{ as required.}\\
		\end{align*}

	\end{proof}

This allows us to apply the enriched small object argument on the sets $\underline{\gencof}$ and $\underline{\gentrivcof}$ of maps in $\CatE$, in which $\gencof$ and $\gentrivcof$ are the sets of maps in $\Cat$ defined in \cref{equ: I and J}.
 
	\begin{cor}
		\label{cor:soa}
		Let $\E$ be a lextensive cartesian closed category with pullback stable coequalisers in which the forgetful functor $\mathcal{U}: \CatE \to \mathbf{Gph}(\E)$ has left adjoint. Then we have weak factorisation systems $(^{\lifts}(\underline{\gencof}^{\elifts}), \underline{\gencof}^{\lifts})$ and $(^{\elifts}(\underline{\gentrivcof}^{\elifts}), \underline{\gentrivcof}^{\lifts})$ on $\CatE$.
	\end{cor}
	
	\begin{proof}
	$\CatE$ is finitely complete and cocomplete by \cite[Proposition 3.3]{Hughes2024Elementary} and \cite[Theorem 5.2]{Hughes2024Colimits}. By \cref{prop: CatE tensored and cotensored}, $\CatE$ is copowered and powered by $\E$. Therefore, by Corollary 7.6.4 of \cite{riehl2014categorical} $\CatE$ is finitely $\E$-bicomplete. Therefore $\CatE$ satisfies the conditions for the enriched small object argument given in Proposition~13.4.2 of \cite{riehl2014categorical}. Finally, $\underline{\gencof}$ and $\underline{\gentrivcof}$ are both finite sets. 
	\end{proof}

\subsection{The classical model structure on $\Cat$}\label{subsec: classical mod struc}

Recall that there is a model structure on $\Cat$ with fibrations given by the isofibrations, cofibrations given by functors which are injective-on-objects and weak equivalences given by the equivalences of categories \cite{joyal2006strong}. We call this the \emph{classical model structure} on $\Cat$. This model structure is monoidal with respect to the cartesian product and also cofibrantly generated by the sets

\begin{align}\label{equ: I and J}
	\gencof &:= \{\emptyset \to \mathbf{1}, \mathbf{1}+\mathbf{1} \to \mathbf{2}, \mathcal{P} \to \mathbf{2}\} \\
	\gentrivcof &:= \{\mathbf{1} \to \mathcal{I}\}.
\end{align}

\subsection{The natural model structure on $\CatE$}\label{subsec: nat mod struc}

We recall a model structure on $\CatE$ which we call the \emph{natural model structure}. This was originally constructed in \cite{everaert2005model}, and it was later noted by Lack that it was an example of what he called a `trivial model structure' on the $2$-category of internal categories, internal functors and internal natural transformations \cite[\S 3.5]{lack2006homotopy}.

We start by defining the fibrations and weak equivalences representably; then we give them an internal description which does not rely on the category $\s.$

\begin{define}
     An internal functor $f: \X \to \Y$ is called an \emph{internal equivalence of categories} if for every $\mathbb{A} \in \CatE_0$, the functor $\Hom_{\Cat}(\mathbb{A}, f): \Hom_{\Cat}(\mathbb{A}, \X) \to \Hom_{\Cat}(\mathbb{A}, \Y)$ is essentially surjective on objects and fully faithful. 

     An internal functor $f: \X \to \Y$ is called an $\emph{internal isofibration}$ if for every $\mathbb{A} \in \CatE_0$, the functor $\Hom_{\Cat}(\mathbb{A}, f): \Hom_{\Cat}(\mathbb{A}, \X) \to \Hom_{\Cat}(\mathbb{A}, \Y)$ is an isofibration in $\Cat$.
 \end{define}

The proof of the following is direct from the representable definitions.

\begin{lem}\label{lem: internal isofib}
    \begin{enumerate}
        \item An internal functor $f: \X \to \Y$ is an internal isofibration iff the map $\iso(\X)_1 \to \iso(\Y)_1 \times_{Y_0} X_0$ is a split epimorphism. 

        \item    An internal functor $f: \X \to \Y$ is an internal equivalence of categories if and only if there exists an internal functor $g: \Y \to \X$ and internal natural isomorphisms $gf \cong 1_{\X}$ and $fg \cong 1_{\Y}$. 
    \end{enumerate}

\end{lem}

Proposition 5.9 of \cite{everaert2005model} defines the cofibrations of this model structure as internal functors which on objects have the left lifting property with respect to the split epimorphisms. In \cref{thm: NatModelStructure}, we include an explicit description of the cofibrations.

We note that a constructive analogue to the $(\mathbf{Inj}, \mathbf{Surj})$ weak factorisation system on $\s$ is achieved replace injective maps by \emph{complemented inclusions}, as is done in \cite{gambino2022constructive} in the setting of simplicial sets. We can formulate this concept more generally in a lextensive category $\E$.

\begin{define}
    A morphism $f: A \to B$ in  $\E$ is called a \emph{complemented inclusion} if there exists an object $C$ of $\E$ such that $f$ is isomorphic to $\iota_A : A \hookrightarrow A + C$.
\end{define}

The following is noted in \cite{gambino_henry_sattler_szumilo_2022}.

 \begin{lem}
		\label{lem: compsplit}
		Let $\E$ be a lextensive category. Then complemented inclusion and split epimorphisms form a weak factorisation system. 
\end{lem}

As such, the cofibrations of the model structure can be described as functors which are a complemented inclusion on objects. We arrive at an explicit description of the model structure that does not rely on representable notions.

 \begin{thm}
 \label{thm: NatModelStructure}
     Let $\E$ be a category that is finitely complete such that $\CatE$ is finitely cocomplete. There is a model structure on $\CatE$ which we denote $(\textbf{Weq},\textbf{Cof}, \textbf{Fib})$ in which 
     \begin{itemize}
         \item $\textbf{Weq}$ is the class of internal equivalences of categories.
         \item $\textbf{Cof}$ is the class of functors which are complemented inclusion on objects.
         \item $\textbf{Fib}$ is the class of internal isofibrations.
     \end{itemize}

    The trivial cofibrations are the split monomorphic equivalences and the trivial fibrations are the split epimorphic equivalences. 

    We call this the \emph{natural model structure on $\CatE$}.
 \end{thm}

 We denote the class of trivial fibrations by $\mathbf{TrivFib}$ and the class of trivial cofibrations by $\mathbf{TrivCof}.$ 
 
 \begin{remark}
    We note that  \cite{Hughes2024Colimits} shows that any lextensive cartesian closed category with pullback stable coequalisers in which the forgetful functor $\mathcal{U}: \CatE \to \mathbf{Gph}(\E)$ has left adjoint satisfies the conditions required for the theorem.
 \end{remark}

We note that for $\E =\s$ and assuming the axiom of choice, we recover the classical model structure on $\Cat$. Hence, this gives us an internal version of this model structure, and without the axiom of choice we obtain a constructive version of the classical model structure on $\Cat$.

\section{Cofibrant Generation}\label{sec: cofib gen}

We show that the $\underline{\gencof}$ and $\underline{\gentrivcof}$ generate the model structure of \cref{thm: NatModelStructure}. 

 \begin{remark}
    We also have internal descriptions of the maps in $\underline{\gencof}$ and $\underline{\gentrivcof}$ which do not rely on their description as ordinary categories. Denote the terminal object of $\CatE$ by $\underline{\mathbf{1}}$. We construct $\underline{\mathbf{2}}$ as the internal category with objects $\underline{\mathbf{1}}+ \underline{\mathbf{1}}$, morphisms $\underline{\mathbf{1}}+\underline{\mathbf{1}} + \underline{\mathbf{1}}$, and in general composable $n$-morphisms by $n$ coproducts of $\underline{\mathbf{1}}$. The source, target, identity and composition maps are induced by the universal property of the coproduct. Using this, we construct $\underline{\mathcal{P}}$ and $\underline{\mathcal{I}}$ as pushouts as in \cref{eq: pushout for P}. This recovers the internal categories required.
 \end{remark}

 First, we set up some notation. Let $f: \X \to \Y$ in $\CatE$. Consider the diagrams:

 \begin{equation*}
     \begin{tikzcd}
         X_1 \arrow[dr, "\Delta", dashed] \arrow[ddr, bend right, "1_{X_1}"'] \arrow[drr, bend left, "1_{X_1}"] & & \\
         & \X(\Par) \arrow[r] \arrow[d]\arrow[dr, phantom, "\lrcorner", very near start] & X_1 \arrow[d, "(d_1{,} d_0)"] \\ 
         & X_1 \arrow[r, "(d_1{,} d_0)"' ] & X_0 \times X_0
     \end{tikzcd}
     \begin{tikzcd}
         \X(\Par) \arrow[dr, "f(\Par)", dashed] \arrow[d] \arrow[r]& X_1\arrow[dr, "f_1"] & \\
         X_1 \arrow[dr, "f_1"'] & \Y(\Par) \arrow[d] \arrow[r] \arrow[dr, phantom, "\lrcorner", very near start] & Y_1  \arrow[d, "(d_1{,} d_0)"]\\
         & Y_1 \arrow[r, "(d_1{,} d_0)"'] & Y_0 \times Y_0
     \end{tikzcd}
 \end{equation*}

 \begin{equation*}
     \begin{tikzcd}
         X_1 \arrow[dr, dashed, "w"] \arrow[drr, bend left, "f_1"] \arrow[ddr, bend right, "d_1 \times d_0", swap] & & \\
         & X_0 \times X_0 \times_{Y_0 \times Y_0} Y_1 \arrow[dr, phantom, "\lrcorner", very near start]  \arrow[r, "\pi_1"] \arrow[d, "\pi_2"'] & Y_1 \arrow[d, "d_1 \times d_0"] \\
         & X_0\times X_0 \arrow[r, "f_0 \times f_0", swap] & Y_0 \times Y_0.
     \end{tikzcd}
 \end{equation*}

 \begin{equation*}
     \begin{tikzcd}
         X_1 \arrow[dr, dashed, "q"] \arrow[drr, bend left, "f_1"] \arrow[ddr, bend right, "\Delta", swap] & & \\
         & \X(\mathcal{P}) \times_{\Y(\mathcal{P})} Y_1 \arrow[dr, phantom, "\lrcorner", very near start]  \arrow[r, "\pi_4"] \arrow[d, "\pi_3"] & Y_1 \arrow[d, "\Delta"] \\
         & \X(\mathcal{P}) \arrow[r, "f(\Par)", swap] & \Y(\mathcal{P}).
     \end{tikzcd}
 \end{equation*}

 \begin{define}
    Let $f: \X \to \Y$ be an internal functor. It is called \emph{full} if the map $w$ is a split epi and \emph{faithful} if $w$ is a monomorphism. It is called \emph{fully faithful} if it is full and faithful.
 \end{define}

 \begin{remark}
 \label{rem: representable}
     Note that if $f$ is fully faithful then $w$ is a monomorphism and a split epimorphism and is therefore an isomorphism. Therefore, an internal fully faithful functor is one such that the following commutative square in $\E$ is a pullback:

     \begin{equation*}
         \begin{tikzcd}
             X_1 \arrow[dr, phantom, "\lrcorner", very near start] \arrow[r, "f_1"] \arrow[d, "(d_1{,} d_0)"'] & Y_1 \arrow[d, "(d_1{,} d_0)"] \\
             X_0 \times X_0 \arrow[r, "f_0 \times f_0"'] & Y_0 \times Y_0
         \end{tikzcd}
     \end{equation*}

     The definitions of faithful and fully faithful are therefore representable ones, and so an internal functor $f: \X \to \Y$ is faithful (resp. fully faithful) iff for all $\mathbb{A} \in \CatE$, $\Hom_{\Cat}(\mathbb{A}, f)$ is a faithful (resp. fully faithful) functor. Split epimorphisms are also a representable notion, so $f: \X \to \Y$ is full iff for all $\mathbb{A} \in \CatE,$ the functor $\Hom_{\Cat}(\mathbb{A}, f)$ is a split epimorphism. Assuming the axiom of choice so that surjections are exactly the split epimorphisms, this is equivalent to $\Hom_{\Cat}(\mathbb{A}, f)$ being full. 
 \end{remark}

 \begin{lem}
 \label{lem: faithful}
     Let $f: \X \to \Y$. The map $q: X_1 \to \X(\Par) \times_{\Y(\Par)} Y_1$ is a split epi if and only if $f$ is faithful. In this case, $X_1 \cong \X(\Par) \times_{\Y(\Par)} Y_1$.
 \end{lem}

 \begin{proof}
     Let $\splt(q)$ denote the splitting of $q$. We must show that $w$ is a monomorphism. By representability of monomorphisms, split epimorphisms, and the evaluation map, it is enough to show that for any $E \in \E$, $\Hom_{\s}(E, w)$ is a monomorphism. 

     Let $E \in \E$ and $\alpha, \beta : E \to X_1$ such that $w\alpha=w\beta$. We show that this implies that $\alpha = \beta.$ By assumption, $w\alpha = w\beta,$ so in particular $(d_1, d_0) \alpha = (d_1, d_0)\beta.$ By \cref{lem: parallel}, we therefore have an induced arrow $(\alpha, \beta): E \to \X(\Par).$ Similarly, $w\alpha=w\beta$,implies that $f_1\alpha = f_1 \beta$, so we have an induced arrow $\Delta f_1 \alpha: E \to \Y(\Par)$ that factors through $Y_1$. It is clear from the universal property of the pullback that $f(\Par)(\alpha, \beta)= \Delta(f_1 \alpha),$ and so we have an induced map $(f_1\alpha , (\alpha, \beta)): E \to \X(\Par) \times_{\Y(\Par)}Y_1.$ The following diagram commutes, exhibiting that $(\alpha, \beta) = \Delta \splt(q) (f_1\alpha, (\alpha, \beta))$, so in particular $\alpha = \splt(q) (f_1\alpha, (\alpha, \beta)) = \beta,$ as required.

     \begin{equation*}
     \begin{tikzcd}
               E \arrow[r, "(f_1\alpha{,} (\alpha{,} \beta))"] \arrow[d, "(\alpha {,} \beta)"] & \X(\Par) \times_{\Y(\Par)} Y_1 \arrow[d, "\pi_3"] \arrow[r, "\text{split}(q)"] & X_1 \arrow[d, "\Delta"] \arrow[r, "q"] & \X(\Par) \times_{\Y(\Par)}Y_1 \arrow[r, "\pi_4"] \arrow[d, "\pi_3"] \arrow[dr, phantom, "\lrcorner", very near start] & Y_1 \arrow[d, "\Delta"] \\
         \X(\Par) \arrow[r, equal] &  \X(\Par) \arrow[r, equal]&  \X(\Par)\arrow[r, equal] & \X(\Par) \arrow[r, "f(\Par)"']& \Y(\Par)
     \end{tikzcd}  
     \end{equation*}

Conversely, suppose $f$ is faithful. We claim that the following diagram commutes providing the desired splitting for $q$. 

\begin{equation*}
\begin{tikzcd}
    \X(\Par) \times_{\Y(\Par)} Y_1 \arrow[r] \arrow[drr, equal] & \X(\Par) \arrow[r, "\pi_1"] & X_1 \arrow[d, "q"] \\
    & & \X(\Par) \times_{\Y(\Par)} Y_1
\end{tikzcd}
\end{equation*}

We use the universal property of $\X(\Par) \times_{\Y(\Par)} Y_1$ as a pullback and show that these arrows agree on the projections to $\X(\Par)$ and $Y_1$. The commutativity of the following diagram witnesses that the arrows agree on the projections to $Y_1$:

\begin{equation*}
    \begin{tikzcd}
        \X(\Par) \times_{\Y(\Par)} Y_1 \arrow[r] \arrow[dd] & \X(\Par) \arrow[r, "\pi_1"] \arrow[d, "f(\Par)"'] & X_1 \arrow[ddr, "f_1"] \arrow[r, "q"] & \X(\Par) \times_{\Y(\Par)}Y_1 \arrow[dd] \\
        & \Y(\Par) \arrow[drr] &  & \\
        Y_1 \arrow[rrr, equal] \arrow[ur, "\Delta"]  & & & Y_1
    \end{tikzcd}
\end{equation*}

To check that they agree on the projection to $\X(\Par),$ we use the universal property of $\X(\Par)$ as a pullback, and show that they agree on the projections $\pi_1: \X(\Par) \to X_1$ and $\pi_2: \X(\Par) \to X_1$. The commutativity of the following diagram shows that the arrows agree on $\pi_1:$

\begin{equation*}
    \begin{tikzcd}
        \X(\Par) \times_{\Y(\Par)} Y_1 \arrow[r] \arrow[d] & \X(\Par) \arrow[r, "\pi_1"]  & X_1 \arrow[ddr, equal] \arrow[r, "q"] & \X(\Par) \times_{\Y(\Par)}Y_1 \arrow[d] \\
       \X(\Par) \arrow[d, "\pi_1"'] & & & \X(\Par) \arrow[d, "\pi_1"] \\
        Y_1 \arrow[rrr, equal]  & & & X_1
    \end{tikzcd}
\end{equation*}

By definition of faithfulness, $w: X_1 \to X_0 \times X_0 \times _{Y_0 \times Y_0} Y_1$ is a monomorphism; to show that the arrows agree on $\pi_2$, we show that they agree on $w_1\pi_2: \X(\Par) \to X_0 \times X_0 \times_{Y_0 \times Y_0} Y_1$ and the result will follow. To show that they agree on this, we use again the universal property of $X_0 \times X_0 \times_{Y_0 \times Y_0} Y_1$ as a pullback and show that they agree on the projections to $X_0 \times X_0$ and $Y_1$. Hence, $q$ is a split epimorphism. It is clearly also a split monomorphism since the following diagram commutes.

    \begin{equation*}
        \begin{tikzcd}
            X_1 \arrow[r] \arrow[drr, equal] & \X(\Par) \times_{\Y(\Par)}Y_1 \arrow[r] & \X(\Par) \arrow[d, "\pi_1"] \\
            & & X_1 
        \end{tikzcd}
    \end{equation*}

    It follows that $q$ is an isomorphism and $X_1 \cong \X(\Par) \times_{\Y(\Par)}Y_1.$
 \end{proof}

\begin{prop}
 \label{prop: characterisation of TrivFib}
		Let $f: \X \to \Y$ be an internal functor in $\CatE$. The following are equivalent:
		
		\begin{enumerate}
			\item The pullback evaluation with all maps in $\underline{\gencof}$ are split epimorphisms.
			\item $f_0: X_0 \to Y_0$ is a split epimorphism on objects and the canonical maps 
			\begin{align*}
				w:X_1 &\to X_0 \times X_0 \times_{Y_0 \times Y_0} Y_1 \\
				q:X_1 &\to X(\mathcal{P})\times_{Y(\mathcal{P})}Y_1 
			\end{align*}
			are split epimorphisms.
			\item $f$ is a split epimorphism on objects and internally fully faithful.
\item $f$ is a a split epimorphism on objects and an adjoint equivalence of categories.
\item $f$ is a trivial fibration.
   
		\end{enumerate}
	\end{prop}
	
\begin{proof}
		By \cref{rem:eval} $(1) \iff (2)$ follows by definition of pullback evaluation and the evaluation map. The equality of (3) and (4) is given in \cite[Remark 8.4]{Hughes2024Elementary}. The equality of (4) and (5) is definitional. By \cref{lem: faithful}, the equality of (2) and (3) follows easily from the definitions of fully faithfulness. 




  \end{proof}

 \begin{prop}
 \label{prop: C TF Ecofibgen}
     The weak factorisation system $(\mathbf{Cof}, \mathbf{TrivFib})$ is $\E$-cofibrantly generated by $\underline{\gencof}$.
 \end{prop}

 \begin{proof}
     By \cref{cor:soa}, $\underline{\gencof}$ generates an $\E$-enriched weak factorisation system. By \cref{prop: characterisation of TrivFib}, the right class of the underlying weak factorisation system coincides with the right class of $(\mathbf{Cof}, \mathbf{TrivFib})$, and hence these are the same weak factorisation system.
 \end{proof}

We have a similar characterisation for the maps in $(\underline{\gentrivcof})^{\elifts}.$
\begin{prop}
\label{prop: characterisation of isofibrations}
Let $f: \X \to \Y$ be an internal functor in $\CatE$. The following are equivalent:

    \begin{enumerate}
        \item $f$ is an internal isofibration.
        \item $\X(\mathcal{I}) \to X_0 \times_{Y_0} \Y(\mathcal{I})$ is a split epimorphism.
        \item The pullback evaluation with $\mathbf{1} \to \mathcal{I}$ is a split epimorphism.
        \item  $f \in (\underline{\gentrivcof})^{\elifts}$ 
        \item For every $\mathbb{A} \in \CatE, \Hom_{\Cat}(\mathbb{A}, f)$ is an isofibration in $\Cat$.
    \end{enumerate}
\end{prop}

\begin{proof}
    (1) and (2) are the same by the definition given in \cite{everaert2005model}; (2) and (3) are equivalent by \cref{lem:yoneda} and definition of the pullback product; (3) and (4) are equivalent by the definition of enriched lifting. (1) and (5) are equivalent as noted in 3.5 of \cite{lack2006homotopy}.
\end{proof}

 \begin{remark}
     By \cref{lem: ev at I is iso}, characterisation (2) of \cref{prop: characterisation of isofibrations} means that if we restrict to $\GpdE$, then this definition of internal isofibration recovers the definition studied in \cite{niefield2019internal}. 
 \end{remark}

 \begin{prop}
 \label{prop: TC F Ecofibgen}
     The weak factorisation system $(\mathbf{TrivCof}, \mathbf{Fib})$ is $\E$-cofibrantly generated by $\underline{\gentrivcof}$.
 \end{prop}

 \begin{proof}
     We apply the same proof as in \cref{prop: C TF Ecofibgen}.
 \end{proof}

 \begin{prop}
     \label{lem: C closed under tensor}
     The classes $\mathbf{Cof}$ and $\mathbf{TrivCof}$ are closed under cartesian products.
 \end{prop}

 \begin{proof}
     
     An internal functor $f: \X \to \Y$ is a cofibration precisely when it is a complemented inclusion on objects, so $f_0: X_0 \to Y_0$ is isomorphic to an inclusion map $X_0 \hookrightarrow X_0 + C$. Hence, for any $\mathbb{A} \in \CatE$, $(\mathbb{A} \times f)_0 = 1_{A_0} \times f_0$ is isomorphic to a map $A_0 \times X_0 \to A_0 \times (X_0 + C)$. By lextensivity, $A_0 \times (X_0 + C)\cong (A_0 \times X_0) + (A_0 \times C)$, and so $\mathbb{A} \times f$ is a cofibration. 

     An internal functor $f: \X \to \Y$ is a trivial cofibration precisely if it is a split monomorphism which an internal adjoint equivalence and a complemented inclusion on objects. Since complemented inclusion on objects functors are closed under cartesian products, it is enough to show that split monomorphisms and internal equivalences are.  Since these and the cartesian product are representable, it is enough to show that this holds in $\Cat$, but this follows from $(\Cat, \times, \mathbf{1})$ being a monoidal model category.
 \end{proof}

 \begin{prop}
 \label{thm: cofibgen}
     The natural model structure on $\CatE$ is cofibrantly generated by the sets $\underline{\gencof}$ and $\underline{\gentrivcof}$. 
 \end{prop}

 \begin{proof}
     By \cref{prop: C TF Ecofibgen,prop: TC F Ecofibgen}, the weak factorisation systems $(\mathbf{Cof}, \mathbf{TrivFib})$ and $(\mathbf{TrivCof}, \mathbf{Fib})$ are $\E$-cofibrantly generated. By \cref{lem: C closed under tensor}, the left classes of these are closed under cartesian product; in particular, they are closed under tensoring with objects in $\E$. By Proposition 13.4.2 of \cite{riehl2014categorical}, this happens if and only if the weak factorisation systems are cofibrantly generated in the ordinary sense by $\underline{\gencof}$ and $\underline{\gentrivcof}$. 
 \end{proof}

 All cofibrantly generated model structures form algebraic model structures, as defined in \cite{riehl2011algebraic}. Therefore, we have the following:

  \begin{cor}
  \label{cor: alg model structure existence}
     There is an algebraic model structure on $\CatE$ with underlying model structure the natural model structure on $\CatE$.
 \end{cor}

We show that this model structure interacts nicely with the cartesian product.  

 \begin{prop}
 \label{thm: sym mon model structure}
     The natural model structure on $(\CatE, \times, \TE)$ is a symmetric monoidal model structure.
 \end{prop}

\begin{proof}
    Since every object is cofibrant and the monoidal product is symmetric, the condition for a monoidal model category simplifies to requiring that for every $\mathbb{A} \in \CatE$, $\mathbb{A} \times -$ preserves cofibrations and trivial cofibrations. This is proven in \cref{lem: C closed under tensor}.
\end{proof}

We arrive at main result of this section. 

\begin{thm}\label{thm: internal analogue}
    Let $\E$ be a lextensive cartesian closed category with pullback stable coequalisers in which the forgetful functor $\mathcal{U}: \CatE \to \mathbf{Gph}(\E)$ has left adjoint. The natural model structure on $\CatE$ has the following properties:

    \begin{enumerate}
        \item It is symmetric monoidal with respect to the cartesian product.
        \item It is cofibrantly generated by the sets
\begin{align*}
	\underline{\gencof} &:= \{\underline{\emptyset} \to \underline{\mathbf{1}}, \underline{\mathbf{1}}+\underline{\mathbf{1}
    }\to \underline{\mathbf{2}}, \underline{\mathcal{P}} \to \underline{\mathbf{2}}\} \\
	\underline{\gentrivcof} &:= \{\underline{\mathbf{1}} \to \underline{\mathcal{I}}\}.
\end{align*}
\item There is an algebraic model structure on $\CatE$ together with the equivalences of categories whose underlying model structure is the natural model structure.

    \end{enumerate}
\end{thm}

\subsection{A model structure on internal modules}

The nice properties of being a cofibrantly generated and symmetric monoidal model structure allow the development of a lot of nice theory, as explored in \cite{hovey1998monoidal}. One such result allows us to lift the model structure on $\CatE$ onto the category of internal modules for an internal monoidal category in $\CatE$.

Recall that a monoid in $\CatE$ is an internal monoidal category. Let $(\mathbb{M}, \otimes, u)$ be an internal monoidal category. An \emph{$\mathbb{M}$-module} is an $\X \in \CatE$ together with an internal functor $\mu: \mathbb{M} \times \X \to \X$ such that the following diagrams commute. 

\begin{equation*}
    \begin{tikzcd}
        \mathbb{M} \times \mathbb{M} \times \X \arrow[r, "\text{id}_{\mathbb{M}} \times \mu"] \arrow[d, "\otimes \times \text{id}_{\X}"'] & \mathbb{M} \times \X \arrow[d, "\mu"] \\
        \mathbb{M} \times \X \arrow[r, "\mu"'] & \X
    \end{tikzcd}\qquad
        \begin{tikzcd}
        \TE \arrow[r, "u \times x"] \arrow[dr, "\forall x"'] & \mathbb{M}\times \X \arrow[d, "\mu"] \\
        & \X
    \end{tikzcd}
\end{equation*}

 A morphism of left $\mathbb{M}$-modules $f: (\X, \mu_{\X}) \to (\Y, \mu_{\Y})$ is an internal functor $f: \X \to \Y$ such that $\mu_{\Y}(1\times f) = f \mu_{\X}. $ This defines a category $\mathbb{M}\text{-}\mathbf{Mod}$. Suppose that $\mathbb{M}$ is furthermore a symmetric monoidal category. Then we can define $\mathbb{M}\text{-}\mathbf{Algebra}$ to be the category of monoids internal to $\mathbb{M}\text{-}\mathbf{Mod}.$ 

\begin{thm}
\label{cor: model structure on modules}
     Let $\E$ be a lextensive cartesian closed category with pullback stable coequalisers in which the forgetful functor $\mathcal{U}: \CatE \to \mathbf{Gph}(\E)$ has left adjoint. Let $\mathbb{M}$ be a monoidal internal category in $\CatE$. Then there is a cofibrantly generated monoidal model structure on the category $\mathbb{M}\text{-}\mathbf{Mod}$ in which a map is a
weak equivalence or fibration if and only if it is a equivalence of internal categories or an internal isofibration. Furthermore, a cofibration of $\mathbb{M}$-modules is a functor which is a complemented inclusion on objects. 

Additionally, if $\mathbb{M}$ is an internal symmetric monoidal category, the category $\mathbb{M}\text{-}\mathbf{Algebra}$ is a model category.
\end{thm}

\begin{proof}
     The first result follows by an application of Proposition 2.8 of \cite{hovey1998monoidal}. Since every object is cofibrant in $\mathbb{M}\text{-}\mathbf{Mod},$ we can apply Theorem 3.3 from \cite{hovey1998monoidal} to get an actual model structure on $\mathbb{M}\text{-}\mathbf{Algebra}.$
\end{proof}

 \section{Algebraic Aspects}\label{sec: Algebraic aspects}

 By \cref{thm: internal analogue} (3), we have an algebraic model structure on $\CatE$. In this section, we provide an explicit description of the two algebraic 
 weak factorisation systems of this, leveraging \cref{thm: internal analogue} (3) to simplify calculations. First, we recall a well-known result from monad theory.

 \begin{define}
     Let $(T: \C \to \C, \eta: \text{id}_{\C} \Rightarrow T)$ be a pointed endofunctor and suppose there exists a left adjoint $F$ to the forgetful functor $G: (T, \eta)\text{-}\mathbf{Alg} \to \C.$ We call $GF$ the \emph{algebraically free monad on} $(T, \eta)$.
 \end{define}

The endofunctor $GF$ does indeed form a monad \cite[Proposition 9.4.1]{barr2000toposes}. For us, the most important property of the algebraically free monad of a pointed endofunctor will be that it has the same algebras as the pointed endofunctor. A proof of this result is given in \cite[Proposition 9.4.5]{barr2000toposes}.
 \begin{prop}
 \label{prop:monadsalgebraicallyfreeoverpointedendofunctor}
     Let $(T, \eta)$ be a pointed endofunctor on a category $\C$, and consider a functor $F:~\C \to (T, \eta)\text{-}\mathbf{Alg}$ given on objects by $X \mapsto (TX, \mu_X:  T^2X \to TX)$. The following are equivalent:

     \begin{enumerate}
         \item $F$ is left adjoint to the forgetful functor $(T, \eta)\text{-}\mathbf{Alg} \to \C$.
         \item $\mathbb{T} : = (T, \eta, \mu)$ is the algebraically free monad over $(T, \eta)$. 
     \end{enumerate}

     In this case, Then $\mathbb{T}\text{-}\mathbf{Alg} \cong (T, \eta)\text{-}\mathbf{Alg}$
 \end{prop}

    \subsection{Factorisations}\label{sec: Functorial Factorisations}

    We now fix the setting of $\C = \CatE^{\mathbf{2}}$. In this subsection, we describe the factorisations of \cite[\S 7]{everaert2005model} which utilise Ken Brown's factorisation lemma \cite{brown1973abstract}.
 Let $f: \X \to \Y$ be an internal functor. We construct the following factorisation

 \begin{equation*}
     \begin{tikzcd}
         \X \arrow[rr, "f"] \arrow[dr, "\mathbf{C}f"'] & & \Y \\
         & \mathbb{E}(f) \arrow[ur, "\mathbf{TF}(f)"'] &
     \end{tikzcd}
 \end{equation*}

 in which $\mathbf{C}(f)$ is a cofibration and $\mathbf{TF}(f)$ is a trivial fibration. We begin on objects:

  \begin{equation*}
     \begin{tikzcd}
         X_0 \arrow[rr, "f"] \arrow[dr, "\mathbf{C}f: = \iota_{X_0}"', hookrightarrow] & & Y_0 \\
         & X_0 + Y_0 \arrow[ur, "\mathbf{TF}(f): = \iota_{X_0} + \text{id}_{Y_0}"'] &
     \end{tikzcd}
 \end{equation*}

    Write $E(f)_0: = X_0 + Y_0$. On morphisms:

    \begin{equation*}
        \begin{tikzcd}[column sep = huge]
            X_1 \arrow[d, "(d_1{,} d_0)"'] \arrow[dr, "\mathbf{C}(f)_1", dashed] \arrow[drr, "f_1", bend left] & & \\ 
            X_0 \times X_0 \arrow[dr, "\mathbf{C}(f)_0 \times \mathbf{C}(f)_0"'] & E(f)_1 \arrow[dr, phantom, "\lrcorner", very near start] \arrow[d, "(d_1{,} d_0)"'] \arrow[r, "\TF(f)_1"] & Y_1 \arrow[d, "(d_1{,}d_0)"] \\
            & E(f)_0 \times E(f)_0 \arrow[r, "\TF(f)_0 \times \TF(f)_0"'] & Y_0 \times Y_0
        \end{tikzcd}
    \end{equation*}

    With this, define the following:

    \begin{equation*}
        \begin{tikzcd}[column sep = huge]
              E(f)_0 \arrow[ddr, "\left(\mathbf{1}_{E(f)_0}{,}\mathbf{1}_{E(f)_0}\right)"', bend right] \arrow[dr, "i", dashed] \arrow[r, "\TF(f)_0"] & Y_0 \arrow[dr, "i"] & \\ 
             & E(f)_1 \arrow[dr, phantom, "\lrcorner", very near start] \arrow[d, "(d_1{,} d_0)"'] \arrow[r, "\TF(f)_1"] & Y_1 \arrow[d, "(d_1{,}d_0)"] \\
            & E(f)_0 \times E(f)_0 \arrow[r, "\TF(f)_0 \times \TF(f)_0"'] & Y_0 \times Y_0
        \end{tikzcd}
    \end{equation*}

    \begin{equation*}
    \begin{tikzcd}
        E(f)_2 \arrow[r] \arrow[dr, phantom, "\lrcorner", very near start] \arrow[d] & E(f)_1 \arrow[d, "d_0"] \\
        E(f)_1 \arrow[r, "d_1"'] & E(f)_0
    \end{tikzcd}
    \end{equation*}

    \begin{equation*}
        \begin{tikzcd}[column sep = huge]
              E(f)_2 \arrow[d] \arrow[dr, "m", dashed] \arrow[r] &E(f)_0 \arrow[dr, "\TF(f)_0"] & \\ 
             E(f)_0 \arrow[dr, "\left(\mathbf{1}_{E(f)_0}{,}\mathbf{1}_{E(f)_0}\right)"'] & E(f)_1 \arrow[dr, phantom, "\lrcorner", very near start] \arrow[d, "(d_1{,} d_0)"'] \arrow[r, "\TF(f)_1"] & Y_1 \arrow[d, "(d_1{,}d_0)"] \\
            & E(f)_0 \times E(f)_0 \arrow[r, "\TF(f)_0 \times \TF(f)_0"'] & Y_0 \times Y_0
        \end{tikzcd}
    \end{equation*}

    With this data, $\mathbb{E}(f) : =  (E(f)_0, E(f)_1, d_0, d_1, i, m)$ forms an internal category and $\mathbf{C}(f) : = (\mathbf{C}(f)_0, \mathbf{C}(f)_1)$ forms an internal functor which is a complimented inclusion on objects by construction and $\TF(f): = (\TF(f)_0, \TF(f)_1)$ forms an internal functor which is split epi on objects and fully faithful by construction. This is the factorisation given for the $(\mathbf{Cof}, \mathbf{TrivFib})$ weak factorisation system.

    The $(\mathbf{TrivCof}, \mathbf{Fib})$ factorisation requires the first factorisation system. For an internal functor $f: \X \to \Y$, we first form the following induced map:

    \begin{equation*}
        \begin{tikzcd}
            \X \arrow[dr, dashed, "\mathbf{W}(f)", near end] \arrow[r] \arrow[ddr, bend right, "(\mathbf{1}{,} f)"'] & \X^{\mathcal{I}} \arrow[dr, "f^{\mathcal{I}}"] & \\
            & \Map(f) \arrow[r] \arrow[d, "p_{\X \times \Y}"'] \arrow[dr, phantom, "\lrcorner", very near start] & \Y^{\mathcal{I}} \arrow[d] \\
            & \X \times \Y \arrow[r, "f \times 1"'] & \Y \times \Y.
        \end{tikzcd}
    \end{equation*}

    Then factorise $\mathbf{W}(f)$ using the $(\mathbf{Cof}, \mathbf{TrivFib})$ weak factorisation system and define $\mathbf{TC}(f) : = \mathbf{C}(\mathbf{W}(f)): \X \to \mathbb{E}(\mathbf{W}(f))$ and $\mathbf{F}(f)$ to be the composite  

    \begin{equation*}
        \begin{tikzcd}[column sep = huge]
            \mathbb{E}(\mathbf{W}(f)) \arrow[r, "\TF(\mathbf{W}(f))"] & \Map(f) \arrow[r, "p_{\X \times \Y}"] & \X \times \Y \arrow[r, "p_{\Y}"] & \Y.
        \end{tikzcd}
    \end{equation*}

    Since projections are fibrations and trivial fibrations are fibrations, and fibrations are closed under composition, it follows that $\mathbf{F}(f)$ is a fibration. Proposition 5.10 of \cite{everaert2005model} shows that $\mathbf{W}(f)$ is an equivalence of categories, in particular satisfying $2$-out-of-$3$, and so since $\TF(\mathbf{W}(f))$ is a weak equivalence, it follows that $\mathbf{TC}(f)$ is too; it is also by construction a complemented inclusion on objects.

    \begin{remark}\label{rem: simplified expressions}
        We note that by extensivity we have:

    \begin{align*}\label{eq: simplified E(f)_1}
        E(f)_1 & : = (E(f)_0 \times E(f)_0) \times_{Y_0 \times Y_0} Y_1 \\
        & = ((X_0+Y_0) \times (X_0 + Y_0)) \times_{Y_0 \times Y_0} Y_1 \\
        & \cong ((X_0 \times X_0) + (X_0 \times Y_0) + (Y_0 \times X_0) + (Y_0 \times Y_0) ) \times_{Y_0 \times Y_0} Y_1 \\
        & \cong (X_0 \times X_0)\times_{Y_0 \times Y_0} Y_1 + X_0 \times_{Y_0} Y_1 + X_0 \times_{Y_0} Y_1 + Y_1. 
        \end{align*}
    \end{remark}

    and also $\Map(f) \cong \X \times_{\Y} \Y^{\mathcal{I}}$.

\subsection{Functorial factorisations}

    We upgrade these factorisations into functorial factorisations in the sense of \cite{grandis2006natural}. Consider a commutative square $(u,v): f \to g$. We supply a functor $\mathbb{E}(u,v) : \mathbb{E}(f) \to \mathbb{E}(g)$ such that $(u, \mathbb{E}(u,v)): \mathbf{C}(f) \to \mathbf{C}(g)$ and $(\mathbb{E}(u,v),v) : \TF(f) \to \TF(g)$ are commutative squares. On objects, the map $E(u,v)_0: E(f)_0 \to E(g)_0$ is given by 

    \begin{equation*}
        \begin{tikzcd}
            & A_0 + B_0 & \\
            A_0 \arrow[ur, hook, bend left, "\mathbf{C}(f)_0"] & & B_0 \arrow[ul, hook', bend right] \\
            X_0 \arrow[u, "u_0"] \arrow[r, hook, "\mathbf{C}(f)_0"] & X_0 + Y_0 \arrow[uu, dashed, "E(u{,}v)_0"] & Y_0 \arrow[l, hook'] \arrow[u, "v_0"] 
        \end{tikzcd}
    \end{equation*}

Note that by construction this clearly provides a commutative square $(u_0, E(u,v)_0): \mathbf{C}(f)_0 \to \mathbf{C}(g)_0$ and a map $(E(u,v)_0 , v_0) : \TF(f)_0 \to \TF(g)_0$ by the universal property of the coproduct.

On morphisms, $E(u,v)_1 : E(f)_1 \to E(g)_1$ is induced by the universal property of $E(g)_1$ as a pullback:

 \begin{equation*}
        \begin{tikzcd}[column sep = huge]
              E(f)_1 \arrow[d, "(d_1 {,} d_0)"'] \arrow[dr, "E(u{,}v)_1", dashed] \arrow[r, "\TF(f)_1"] & Y_1 \arrow[dr, "v_1"] & \\ 
             E(f)_0 \times E(f)_0 \arrow[dr, "E(u{,}v)_0 \times E(u{,}v)_0"'] & E(g)_1 \arrow[dr, phantom, "\lrcorner", very near start] \arrow[d, "(d_1{,} d_0)"'] \arrow[r, "\TF(g)_1"] & B_1 \arrow[d, "(d_1{,}d_0)"] \\
            & E(g)_0 \times E(g)_0 \arrow[r, "\TF(g)_0 \times \TF(g)_0"'] & B_0 \times B_0.
        \end{tikzcd}
    \end{equation*}

    Again, by construction this clearly gives a commutative square $(E(u,v)_1, v_1) : \TF(f)_1 \to \TF(g)_1$, and by the universal property of the pullback, it is easy to show that it also gives a commutative square $(u_1, E(u,v)_1): \mathbf{C}(f)_1 \to \mathbf{C}(g)_1$. This assembles into a functor as required; it preserves identities by construction and preservation of composition can be proven using the universal properties involved.

    Hence we have an endofunctor $\mathbf{C}: \CatE^{\mathbf{2}} \to \CatE^{\mathbf{2}}$ with copoint $\epsilon:  \mathbf{C} \Rightarrow \text{id}$ given on $f: \X \to \Y$ by

    \begin{equation*}
        \begin{tikzcd}
            \X \arrow[r, equal] \arrow[d, "\mathbf{C}(f)"'] & \X \arrow[d, "f"] \\
            \mathbb{E}(f) \arrow[r, "\TF(f)"'] & \Y.
        \end{tikzcd}
    \end{equation*}

    Furthermore, we have an endofunctor $\TF: \CatE^{\mathbf{2}} \to \CatE^{\mathbf{2}}$ with point $\eta: \text{id} \Rightarrow \TF$ given on $f: \X \to \Y$ by

    \begin{equation*}
        \begin{tikzcd}
            \X \arrow[r, "\mathbf{C}(f)"] \arrow[d, "f"'] & \mathbb{E}(f) \arrow[d, "\TF(f)"] \\
            \Y \arrow[r, equal] & \Y.
        \end{tikzcd}
    \end{equation*}

    \begin{prop}
        The (co)pointed endofunctors $(\mathbf{C}, \epsilon)$ and $(\mathbf{TF}, \eta)$ provide a functorial factorisation for the weak factorisation system $(\mathbf{Cof}, \mathbf{TrivFib}).$ 
    \end{prop}

    We do the same for the $(\mathbf{TrivCof}, \mathbf{Fib})$ weak factorisation system. It suffices to show that the assignment $f \mapsto \mathbf{W}(f)$ is functorial as then we can define $\mathbf{C} : = \mathbf{TC} \circ \mathbf{W}$ and $\mathbf{F}$ similarly as the composition of functors. Consider the commutative square $(u,v): f \to g.$ We supply a functor $\Map(u,v) : \Map(f) \to \Map(g)$ induced by the universal property of the mapping path space:

    \begin{equation*}
        \begin{tikzcd}
            \Map(f) \arrow[dr, dashed, "\Map(u{,}v)", near end] \arrow[r] \arrow[d] & \Y^{\mathcal{I}} \arrow[dr, "v^{\mathcal{I}}"] & \\
           \X \times \Y \arrow[dr, "u \times v"'] & \Map(g) \arrow[r] \arrow[d] \arrow[dr, phantom, "\lrcorner", very near start] & \B^{\mathcal{I}} \arrow[d] \\
            & \A \times \B \arrow[r, "g \times 1"'] & \B \times \B.
        \end{tikzcd}
    \end{equation*}

This provides a commutative square $(u, \Map(u,v)): \mathbf{W}(f) \to \mathbf{W}(g).$ Respect for identities follows easily from the definition and respect for composition follows from the uniqueness given by the universal property of the pullback. Hence, this assignment is functorial, and we obtain an endofunctor $\mathbf{TC}: \CatE^{\mathbf{2}} \to \CatE^{\mathbf{2}}$ with copoint $\sigma: \mathbf{TC} \Rightarrow \text{id}$ given on $f: \X \to \Y$ by:

  \begin{equation*}
        \begin{tikzcd}
            \X \arrow[r, equal] \arrow[d, "\mathbf{TC}(f)"'] & \X \arrow[d, "f"] \\
            \mathbb{E}(\mathbf{W}(f)) \arrow[r, "\mathbf{F}(f)"'] & \Y.
        \end{tikzcd}
    \end{equation*}

    Similarly, we have an endofunctor $\F: \CatE^{\mathbf{2}} \to \CatE^{\mathbf{2}}$ with point $\xi: \text{id} \Rightarrow \mathbf{F}$ given on $f: \X \to \Y$ by:

    \begin{equation*}
        \begin{tikzcd}
            \X \arrow[r, "\mathbf{TC}(f)"] \arrow[d, "f"'] & \mathbb{E}(\mathbf{W}(f)) \arrow[d, "\mathbf{F}(f)"] \\
            \Y \arrow[r, equal] & \Y.
        \end{tikzcd}
    \end{equation*}

This proves that the weak factorisation system $(\mathbf{TC}, \mathbf{F})$ is functorial. Note that the factorisation of a commutative square $(u,v): f \to g$ is given explicitly as:

\begin{equation*}
    \begin{tikzcd}[column sep = huge]
        \X \arrow[r, "u"] \arrow[d, "\mathbf{TC}(f)"'] & \A \arrow[d, "\mathbf{TC}(g)"] \\
        \mathbb{E}(\mathbf{W}(g)) \arrow[r, "\mathbb{E}(u{,} \Map(u{,}v))"] \arrow[d, "\mathbf{F}(f)"']& \mathbb{E}(\mathbf{W}(f)) \arrow[d, "\mathbf{F}(g)"]\\ 
        \Y \arrow[r, "v"] & \B
    \end{tikzcd}
\end{equation*}

\begin{prop}
    The (co)pointed endofunctors $(\mathbf{TC}, \sigma)$ and $(\mathbf{F}, \xi)$ provide a functorial factorisation for the weak factorisation system $(\mathbf{Cof}, \mathbf{TrivFib}).$ 
\end{prop}

    \subsection{Algebraic structure}\label{subsec: alg structure}

     We unpack the definition of an algebra for a pointed endofunctor to notice a correspondence between algebra structure and fibrational structure. Similarly, there is a correspondence between coalgebras for a copointed endofunctor and cofibrational structure.
     Let $(\Lagr, \R)$ be a functorial weak factorisation system on $\CatE$ with endofunctors $L, R: \CatE^{\mathbf{2}} \to \CatE^{\mathbf{2}}$ with point $\eta: \text{id} \Rightarrow R$ given by 

     \begin{equation*}
         \begin{tikzcd}
             \X \arrow[d, "f"'] \arrow[r, "L(f)"] & \bullet \arrow[d, "R(f)"] \\
             \Y \arrow[r, equal] & \Y
         \end{tikzcd}
     \end{equation*}

     and copoint $\epsilon : L \Rightarrow \text{id}$ given by 

     \begin{equation*}
      \begin{tikzcd}
             \X \arrow[d, "L(f)"'] \arrow[r, equal] & \X \arrow[d, "f"] \\
             \bullet \arrow[r, "R(f)"'] & \Y
         \end{tikzcd}
     \end{equation*}

     An algebra for $(R, \eta)$ is therefore a pair $(f, \overrightarrow{\phi}: R(f) \to f)$ satisfying the following commutativity condition:

     \begin{equation*}
          \begin{tikzcd}
             \X \arrow[d, "f"'] \arrow[r, "L(f)"] & \bullet \arrow[d, "R(f)"] \arrow[r, "\phi"] & \X \arrow[d, "f"] \\
             \Y \arrow[r, equal] & \Y \arrow[r, equal] & \Y
         \end{tikzcd}
     \end{equation*}

     This data can be rearranged into the following commutative diagram:

     \begin{equation*}
         \begin{tikzcd}
             \X \arrow[d, "L(f)"'] \arrow[r, equal] & \X \arrow[d, "f"] \\
             \bullet \arrow[r, "R(f)"'] \arrow[ur, "\phi"] & \Y. 
         \end{tikzcd}
     \end{equation*}

     In other words, $\phi$ provides a lifting of $f$ against its left factor. Conversely, if $f$ lifts against its left factor, then this provides an $(R, \eta)$-algebra structure. 

     The dual of this is also true: $f$ lifts against its right factor if and only if there exists an $\Lagr$-coalgebra structure upon it. It is a theorem due to Garner that for a functorial weak factorisation system $(\Lagr, \R)$, $(R, \eta)$-algebras compose with $(R, \eta)$-algebras coherently if and only if $(L, \epsilon)$ extends to a comonad; dually, $(L, \epsilon)$-coalgebras compose coherently with $(L, \epsilon)$-coalgebras if and only if $(R, \eta)$ extends to a monad \cite[2.24]{riehl2011algebraic}.

     In our setting, we have endofunctors $\mathbf{F}, \TF, \mathbf{C}, \TC$ defined in \cref{sec: Functorial Factorisations}. By the above, these are all (co)pointed; we denote their points by $\xi: \text{id} \Rightarrow \mathbf{F}$, $\eta: \text{id} \Rightarrow \TF$, $\epsilon: \mathbf{C} \Rightarrow \text{id}$ and $\sigma: \TC \Rightarrow \text{id}$ respectively. Therefore, we have the following.

\begin{prop}
\label{prop: algebras for pointed endofunctors}
    Let $f: \X \to \Y$ in $\CatE$.  

    \begin{enumerate}
        \item There exists an $(\mathbf{F}, \xi)$-algebra structure for $f$ is and only if $f$ is a fibration.
        \item There exists an $(\TF, \eta)$-algebra structure for $f$ is and only if $f$ is a trivial fibration.
        \item There exists an $(\mathbf{C}, \epsilon)$-algebra structure for $f$ is and only if $f$ is a cofibration.
        \item There exists an $(\mathbf{TC}, \sigma)$-algebra structure for $f$ is and only if $f$ is a trivial cofibration.
    \end{enumerate}
\end{prop}

    We claim that the (co)pointed endofunctors $(\F, \xi)$ and $(\TF, \eta)$ (resp. $(\mathbf{C}, \epsilon)$ and $(\TC, \sigma)$)  extend to (co)monads; to simplify calculations, we appeal to \cref{cor: alg model structure existence}, which tells us that an algebraic model structure exists, and therefore algebraic weak factorisation systems with these underlying maps exist. 

    For $f: \X \to \Y$, we construct a commutative square $\mu^f: \mathbb{E}(\TF(f)) \to \mathbb{E}(f)$. On objects, it is given by:
    
    \begin{equation*}
        \begin{tikzcd}
            & X_0 + Y_0 & \\
            X_0+ Y_0 \arrow[ur, equal, bend left] \arrow[r, hook] & X_0 + Y_0 + Y_0 \arrow[u, dashed, "\mu^f_0"] & Y_0 \arrow[l, hook'] \arrow[ul, hook', bend right] 
        \end{tikzcd}
    \end{equation*}

    Note that this has the property that $\TF(f)_0 \mu^f_0 = \TF^2(f)_0,$ due to the definition of $\mathbf{TF}^2(f)_0$ as the unique arrow which makes the outside of the following diagram commute:

    \begin{equation}\label{TF^2(f)}
        \begin{tikzcd}
           & Y_0 & \\
            & X_0 + Y_0 \arrow[u, "\TF(f)_0"] & \\
            X_0+ Y_0 \arrow[uur, "\TF(f)_0", bend left] \arrow[ur, equal, bend left] \arrow[r, hook] & X_0 + Y_0 + Y_0 \arrow[u, dashed, "\mu^f_0"] & Y_0 \arrow[l, hook'] \arrow[ul, hook', bend right] \arrow[uul, equal, bend right].
        \end{tikzcd}
    \end{equation}
    
    On morphisms it is given by:

    \begin{equation*}
        \begin{tikzcd}[column sep = huge]
            \mathbb{E}(\TF(f))_1 \arrow[dr, dashed, "\mu^f_1"] \arrow[rrd, "\TF^2(f)_1", bend left] \arrow[d, "(d_1{,} d_0)"] & & \\
            \mathbb{E}(\TF(f))_0 \times \mathbb{E}(\TF(f))_0 \arrow[dr, "\mu^f_0 \times \mu^f_0"'] & \mathbb{E}(f)_1 \arrow[r, "\TF(f)_1"] \arrow[d, "(d_1 {,} d_0)"'] \arrow[dr, phantom, "\lrcorner", very near start] & Y_1 \arrow[d, "(d_1 {,} d_0)"] \\
            & \mathbb{E}(f)_0 \times \mathbb{E}(f)_0 \arrow[r, "\TF(f)_0 \times \TF(f)_0"'] & Y_0 \times Y_0
        \end{tikzcd}
    \end{equation*}

    for which the outer diagram commutes by functoriality of $\TF$ and \cref{TF^2(f)}.

    The map $\mu^f : = (\mu^f_0 , \mu^f_1)$ assembles into an internal functor and hence provides a commutative square $\overrightarrow{\mu_f}: \TF^2(f) \to \TF(f)$ given by 

    \begin{equation*}
        \begin{tikzcd}
            \mathbb{E}(\TF(f)) \arrow[r, "\mu^f"] \arrow[d, "\TF^2(f)"'] & \mathbb{E}(f) \arrow[d, "\TF(f)"] \\
            \Y \arrow[r, equal] & \Y.
        \end{tikzcd}
    \end{equation*}

    Consequently, this forms a natural transformation $\overrightarrow{\mu}: \TF^2\Rightarrow \TF$.

    The following shows that we can give algebraic structure to $\TF(f)$, which exhibits it as a trivial fibration. Through the algebraic lens, this is the same as giving the free $(\TF, \eta)$-algebra structure to $f$. 

    \begin{lem}
        For any $f: \X \to \Y$, $(\TF(f), \overrightarrow{\mu_f})$ is a $(\TF, \eta)$-algebra. 
    \end{lem}

    \begin{proof}
        By construction, it is clear that $\TF(f) \mu^f = \TF^2(f).$ It remains to prove that $\mu^f \Cof(f) = \mathbf{1}_{\mathbb{E}(f)}$. By construction, this is true on objects. On morphisms, we appeal to the universal property of $\mathbb{E}(f)_1$ as a pullback. That is, we show that these both agree on $\TF(f)_1$ and $(d_1, d_0).$ This is witnessed by the commutativity of the following diagrams:

        \begin{equation*}
            \begin{tikzcd}
                \mathbb{E}(f)_1 \arrow[d, equal] \arrow[drr, "\TF(f)_1"'] \arrow[r, "\Cof \TF (f)_1"] & \mathbb{E}(\TF(f))_1 \arrow[r, "\mu^f_1"] \arrow[dr, "\TF^2(f)_1"] & \mathbb{E}(f)_1 \arrow[d, "\TF(f)_1"] \\
                \mathbb{E}(f)_1 \arrow[rr, "\TF(f)_1"']  &  & Y_1
            \end{tikzcd}
        \end{equation*}

        \begin{equation*}
            \begin{tikzcd}
                \mathbb{E}(f)_1 \arrow[dd, equal]  \arrow[dr, "(d_1 {,} d_0)"] \arrow[rr, "\Cof \TF (f)_1"] & & \mathbb{E}(\TF(f))_1 \arrow[d, "(d_1 {,} d_0)"] \arrow[r, "\mu^f_1"]& \mathbb{E}(f)_1 \arrow[dd, "(d_1{,} d_0)"] \\
                & \mathbb{E}(f)_0 \times \mathbb{E}(f)_0 \arrow[r, hook, "\Cof \TF (f)_0"] \arrow[drr, equal]& \mathbb{E}(\TF(f))_0 \times \mathbb{E}(\TF(f))_0 \arrow[dr, "\mu^f_0 \times \mu^f_0"] & \\
                \mathbb{E}(f)_1 \arrow[rrr, "(d_1 {,} d_0)"']  & & & \mathbb{E}(f)_0 \times \mathbb{E}(f)_0 
            \end{tikzcd}
        \end{equation*}
    \end{proof}

    Consider the functor $\overline{\TF}: \CatE^{\mathbf{2}} \to (\TF, \eta)\text{-}\mathbf{Alg}$ defined by $f \mapsto (\TF(f), \overrightarrow{\mu_f})$. To show that this is well-defined on morphisms, one must show that for $(u,v): f \to g$, $\mathbb{E}(u,v)\mu^f = \mu^g \mathbb{E}(\mathbb{E}(u,v),v)$, which follows from the universal properties of $E(\TF(f))_0$ as a coproduct to show equality on objects and of $E(g)_1$ as a pullback to show equality on morphisms.

    \begin{prop}
        $\overline{\TF}$ is left adjoint to the forgetful functor $(\TF, \eta)\text{-}\mathbf{Alg} \to \CatE^{\mathbf{2}}.$ Consequently, $\mathbb{TF} := (\TF, \eta, \overrightarrow{\mu})$ is a monad. 
    \end{prop}

    \begin{proof}
        We describe mutually inverse functors 
         \begin{equation*}
        \begin{tikzcd}
            \Theta : (\TF, \eta)\text{-}\mathbf{Alg}((Tf, \overrightarrow{\mu_f}), (g, \alpha)) \arrow[r, shift right] & \CatE^{\mathbf{2}}(f,U(g, \alpha)) \arrow[l, shift right] : \Lambda.
        \end{tikzcd}
    \end{equation*}

    Given a morphism $\phi: (\TF(f), \overrightarrow{\mu_f}) \to (g, \alpha),$ we define $\Theta(\phi): f \to g$ by the composition 

    \begin{equation*}
        \begin{tikzcd}
            f \arrow[r, "\eta_f"]  & \TF(f) \arrow[r, "\phi"] & g
        \end{tikzcd}
    \end{equation*}

    Conversely, given $\psi: f \to U(g, \alpha)$, we give a homomorphism of algebras $\Lambda(\psi) : (\TF(f), \overrightarrow{\mu_f}) \to (g, \alpha)$ by the composite

    \begin{equation*}
        \begin{tikzcd}
            \TF(f) \arrow[r, "\TF(\psi)"]  & \TF(g) \arrow[r, "\alpha"] & g
        \end{tikzcd}
    \end{equation*}

    To show this is a well-defined homomorphism of algebras, one must show that $$\alpha\TF(\alpha\TF(\psi)) = \alpha\TF(\psi)\overrightarrow{\mu_f},$$ which can be done using the universal properties of $\mathbb{E}(\TF(f))$ and $\mathbb{E}(g)$. To show that these maps are mutually inverse, one leverages the same universal properties and the unit property for $\alpha$ to be the structure of an algebra. We omit the tedious details.
    \end{proof}

    This process, or the dual of it, works to extend the (co)pointed endofunctors $\Cof, \TC$ and $\F$ into (co)monads. Rather than repeat the explanations, we just give the construction of the (co)algebraic structure and leave the details to the interested reader.

    For $f: \X \to \Y$, we construct an internal functor $\delta^f : \mathbb{E}(f) \to \mathbb{E}(\Cof(f)).$ On objects it is given by :

    \begin{equation*}
        \begin{tikzcd}
            & X_0 + X_0 + Y_0 & \\
            X_0 \arrow[r, hook] \arrow[ur, hook', bend left ] & X_0 + Y_0 \arrow[u, dashed, "\delta^f_0"] & Y_0 \arrow[l, hook'] \arrow[ul, hook, bend right] 
        \end{tikzcd}
    \end{equation*}

    On morphisms, it is given by: 

        \begin{equation*}
        \begin{tikzcd}[column sep = huge]
            \mathbb{E}(f)_1 \arrow[dr, dashed, "\delta^f_1"] \arrow[rrrd, equal, bend left] \arrow[d, "(d_1{,} d_0)"] &&  & \\
            \mathbb{E}(f)_0 \times \mathbb{E}(f)_0 \arrow[dr, "\delta^f_0 \times \delta^f_0"'] & \mathbb{E}(\Cof(f))_1 \arrow[rr, "\TF(\Cof(f))_1"] \arrow[d, "(d_1 {,} d_0)"'] \arrow[drr, phantom, "\lrcorner", very near start] & &  \mathbb{E}(f)_1 \arrow[d, "(d_1 {,} d_0)"] \\
            & \mathbb{E}(\Cof(f))_0 \times \mathbb{E}(\Cof(f))_0 \arrow[rr, "\TF(\Cof(f))_0 \times \TF(\Cof(f))_0"'] & & \mathbb{E}(f)_0 \times \mathbb{E}(f)_0
        \end{tikzcd}
    \end{equation*}

    We obtain a morphism $\overrightarrow{\delta_f}: \Cof(f) \to \Cof^2(f)$ given by

    \begin{equation*}
        \begin{tikzcd}
            \X \arrow[r, equal] \arrow[d, "\Cof(f)"'] & \X \arrow[d, "\Cof^2(f)"] \\
            \mathbb{E}(f) \arrow[r, "\delta^f"'] & \mathbb{E}(\Cof(f)).
        \end{tikzcd}
    \end{equation*}

    We therefore have a natural transformation $\overrightarrow{\delta}: \Cof \Rightarrow \Cof^2$.

    \begin{lem}
        For any $f: \X \to \Y$, $(f, \overrightarrow{\delta_f})$ is a $(\Cof, \epsilon)$-coalgebra.
    \end{lem}

     Consider the functor $\overline{\Cof}: \CatE^{\mathbf{2}} \to (\Cof, \epsilon)\text{-}\mathbf{Coalg}$ defined by $f \mapsto (\Cof(f), \overrightarrow{\delta_f}).$ 

    \begin{prop}
        $\overline{\Cof}$ is right adjoint to the forgetful functor $(\Cof, \epsilon)\text{-}\mathbf{Coalg} \to \CatE^{\mathbf{2}}.$ Consequently, $\mathbb{C} : = (\Cof, \epsilon, \overrightarrow{\delta})$ is a comonad. 
    \end{prop}

    The constructions involved in the analogous procedure for $(\F, \xi)$ and $(\TC, \sigma)$ are a little bit more involved. We start by treating $(\F, \xi).$ Firstly, we define a map $\beta~:~E(\Weq(f))_0 \to X_0$ given by 

    \begin{equation*}
        \begin{tikzcd}
            & X_0 & \\
            X_0 \arrow[r, hook] \arrow[ur, equal] & X_0 + \Map(f)_0 \arrow[u, dashed, "\beta"] & \Map(f)_0 \arrow[l, hook'] \arrow[ul, "\pi_{X_0}"']. 
        \end{tikzcd}
    \end{equation*}

    Note that $f_0 \beta = \mathbf{F}(f)_0$ by the universal property of $E(\Weq(f))_0$ as a coproduct. 
    
    We use this to define a map $\tau: \Map(\F(f))_0 \to \Map(f)_0$:

    \begin{equation*}
        \begin{tikzcd}
            \Map(\F(f))_0 \arrow[d] \arrow[dr, dashed, "\tau"] \arrow[r] & \iso(\Y)_1 \arrow[dr, equal] & \\
            E(\Weq(f))_0 \times Y_0 \arrow[dr, "\beta \times \mathbf{1}"'] & \Map(f)_0 \arrow[r] \arrow[d] \arrow[dr, phantom, "\lrcorner", very near start] & \iso(\Y)_1 \arrow[d] \\
            & X_0 \times Y_0 \arrow[r, "f \times \mathbf{1}"'] & Y_0 \times Y_0
        \end{tikzcd}
    \end{equation*}

    for which the outer diagram commute by the commutativity of the following:

      \begin{equation*}
        \begin{tikzcd}
            \Map(\F(f))_0 \arrow[d]  \arrow[r] & \iso(\Y))_1 \arrow[dr, equal] \arrow[d] & \\
            E(\Weq(f))_0 \times Y_0 \arrow[dr, "\beta \times \mathbf{1}"'] \arrow[r, "\F(f)_0 \times \mathbf{1}"] & Y_0 \times Y_0 \arrow[dr, equal] & \iso(\Y)_1 \arrow[d] \\
            & X_0 \times Y_0 \arrow[r, "f_0 \times \mathbf{1}"'] & Y_0 \times Y_0.
        \end{tikzcd}
    \end{equation*}

    Next, we construct a map $\kappa^f: \mathbb{E}(\Weq\F(f)) \to \mathbb{E}(\Weq(f)).$ On objects:

    \begin{equation*}
        \begin{tikzcd}[row sep= small]
            & X_0 + \Map(f)_0 & \\
            & & \Map(f)_0 \arrow[ul, hook] \\
            X_0 + \Map(f)_0 \arrow[uur, equal] \arrow[r, hook] & X_0 + \Map(f)_0 + \Map(\F(f))_0 \arrow[uu, dashed, "\kappa^f_0"] & \Map(\F(f))_0 \arrow[l, hook'] \arrow[u, "\tau"']
         \end{tikzcd}
    \end{equation*}

    Note that $\F(f)_0\kappa^f_0 = \F^2(f)_0$, which can be check using the universal property of $E(\Weq\F(f))_0$ as a coproduct.

    We extend this to morphisms by the universal property of $E(\Weq(f))_1$ as a pullback:

    \begin{equation*}
        \begin{tikzcd}[column sep = huge]
            E(\Weq\F(f))_1 \arrow[rrd, "\F^2(f)_1", bend left] \arrow[d, "(d_1 {,} d_0)"']  \arrow[dr, dashed, "\kappa^f_1"]& & \\
            E(\Weq\F(f))_0 \times E(\Weq \F (f))_0 \arrow[dr, "\kappa^f_0 \times \kappa^f_0"'] & E(\Weq(f))_1 \arrow[r, "\F(f)_1"] \arrow[d, "(d_1 {,} d_0)"'] \arrow[dr, phantom, "\lrcorner", very near start] & Y_1 \arrow[d, "(d_1{,}d_0)"] \\
           & E(\Weq(f))_0 \times E(\Weq(f))_0 \arrow[r, "\F(f)_0 \times \F(f)_0"'] & Y_0 \times Y_0
        \end{tikzcd}
    \end{equation*}

    This provides a commutative square $\overrightarrow{\kappa_f}: \F^2(f) \to \F(f)$ given by

    \begin{equation*}
        \begin{tikzcd}
            \mathbb{E}(\Weq\F(f)) \arrow[d, "\F^2(f)"'] \arrow[r, "\kappa^f"] & \mathbb{E}(\Weq(f)) \arrow[d, "\F(f)"] \\
            \Y \arrow[r, equal] & \Y. 
        \end{tikzcd}
    \end{equation*}

    and therefore a natural transformation $\overrightarrow{\kappa}: \F^2 \Rightarrow \F$.

    \begin{lem}
        For any $f: \X \to \Y$, $(\F(f), \overrightarrow{\kappa_f})$ is an $(\F, \xi)$-algebra.
    \end{lem}

    Consider the functor $\overline{\F}: \CatE^{\mathbf{2}} \to (\F, \xi)\text{-}\mathbf{Alg}$ defined by $f \mapsto (\F(f), \overrightarrow{\kappa_f})$.

    \begin{prop}
    \label{prop: F monad}
        $\overline{\F}$ is left adjoint to the forgetful functor $(\F, \xi)\text{-}\mathbf{Alg} \to \CatE^{\mathbf{2}}.$ Consequently, $\mathbb{F}: = (\F, \xi, \overrightarrow{\kappa})$ is a monad. 
    \end{prop}

    Whilst it is possibly to do the same for $(\mathbf{TC}, \sigma)$ and construct algebraic structure, it requires many lengthy and unenlightening calculations. Moreover, we will not actually need the explicit description of this comonad, so we instead appeal to a universal property and prove that $(\mathbf{TC}, \sigma)$ extends to a comonad by showing that it is the underlying left factor of an algebraic weak factorisation system. 

    \begin{prop}
        The monads $\mathbb{F}$ and $\mathbb{TF}$ are the right classes of algebraic weak factorisation systems.
    \end{prop}

    \begin{proof}
         We shall treat $\mathbb{F}$; the other case is similar. We apply \cref{thm: internal analogue} (3); since the weak factorisation system $(\mathbf{TrivCof}, \mathbf{Fib})$ is cofibrantly generated, we can upgrade it to an algebraic weak factorisation system by Garner's algebraic version of the small object argument; this is $(\mathbb{L}, \mathbb{R})$ such that $\mathbb{R}$ is a monad whose $\mathbb{R}$-algebras are exactly the fibrations. This monad has the universal property of being algebraically free on a pointed endofunctor $\mathbf{R}$ whose algebras are fibrations. By \cref{prop: algebras for pointed endofunctors,prop: F monad}, this is exactly what $\mathbb{F}$ is, whence $\mathbb{F} \cong \mathbb{R}$ and so $\mathbb{F}$ is the right class of an algebraic weak factorisation system. 
    \end{proof}

    \begin{cor}
        The algebraically cofree comonad on the copointed endofunctor $(\mathbf{TC}, \sigma)$ exists.
    \end{cor}

    \begin{proof}
        $\mathbb{F}$ is the right class of some algebraic weak factorisation system $(\mathbb{T}, \mathbb{F})$, where $\mathbb{T}=(\mathbf{T}, a, b)$. This algebraic weak factorisation system has underlying ordinary weak factorisation system $(\mathbf{T}, \mathbf{F}).$ Ordinary weak factorisation systems are determined by their left class, so it follows that $\mathbf{T} = \mathbf{TC}$ and $a = \sigma,$ as required.  
    \end{proof}

We denote the algebraically cofree comonad on $(\mathbf{TC}, \sigma)$ by $\mathbb{TC} : = (\mathbf{TC}, \sigma, \rho).$ Note that we could have used the same reasoning to construct $\mathbb{C},$ but in this case it was simple to do directly. 

\begin{cor}
    \label{prop: awfs}
       Both $(\mathbb{TC}, \mathbb{F})$ and $(\mathbb{C}, \mathbb{TF})$ are algebraic weak factorisation systems. 
\end{cor}

    \begin{thm}
    \label{thm: algebraic model structure on CatE}
        Let $\E$ be a lextensive cartesian closed category with pullback stable coequalisers in which the forgetful functor $\mathcal{U}: \CatE \to \mathbf{Gph}(\E)$ has left adjoint. The algebraic weak factorisation systems $(\mathbb{TC}, \mathbb{F})$ and $(\mathbb{C}, \mathbb{TF})$ provide an algebraic model structure on $(\CatE, \mathbf{Weq})$ whose underlying model structure is the natural model structure on $\CatE$. 
    \end{thm}

    \begin{proof}
        \cref{prop: awfs} tells us that these algebraic weak factorisation systems are isomorphic as algebraic weak factorisation systems to the ones created from Garner's small object argument. The result then follows from \cref{cor: alg model structure existence} and \cref{thm: NatModelStructure}.
    \end{proof}

    An important result for us is that this algebraic model structure restricts to the category of internal groupoids. We prove this below.

    \begin{thm}
    \label{cor: algebraic model structure on GpdE}
        Let $\E$ be a lextensive cartesian closed category with pullback stable coequalisers in which the forgetful functor $\mathcal{U}: \CatE \to \mathbf{Gph}(\E)$ has left adjoint. The algebraic natural model structure on $\CatE$ restricts to an algebraic natural model structure on $\GpdE$. 
    \end{thm}

    \begin{proof}
        We need to show that given $f: \mathbb{G} \to \mathbb{H}$ in $\GpdE$, the intermediate objects $\mathbb{E}(f)$ and $\mathbb{E}(\Weq(f))$ are in $\GpdE$. We construct $E(f)_1$ by pulling back over $H_1 = \iso(\mathbb{H})_1$.

        Similarly, we construct $E(\Weq(f))_1$ by pulling back over $\Map(f)_1$, which contains only isomorphisms as $\Map(f)$ is constructed as a pullback and $\GpdE$ and $\iso$ is a right adjoint so preserves limits, yielding $\iso(\Map(f))_1 = \Map(f)_1$, as required. 
    \end{proof}

\subsection{Equivalent algebraic descriptions}

In the previous section, we described the (co)algebras for (co)monads that form an algebraic model structure. In this section, we give equivalent descriptions of these algebras that are more useful and familiar. We show that the structure of being an $\mathbb{F}$-algebra is equivalent to the structure of being a cloven isofibration, and the structure of being a $\mathbb{TC}$-coalgebra is equivalent to the structure of being a retract complemented inclusion on objects with the structure of a $2$-cell witnessing the fact that is an equivalence. It is also true that the structure of being a $\mathbb{C}$-coalgebra is equivalent the structure of being a complemented inclusion on objects and the structure of being a $\mathbb{TF}$-algebra is equivalent to the structure of a splitting on objects and a $2$-cell witnessing the fact that it is an equivalence.

We first treat $\mathbb{TC}$-coalgebras.

\begin{define}
\label{def: algret}
    Let $g:\A \to \Y$ be an internal functor, $r: \Y \to \A$ be a retract of $g$, $j: X_0 + C \to Y_0$ be an isomorphism such that $j \cdot \iota_{X_0} = f_0$, and $\beta: gr \Rightarrow \text{id}_{\Y}$ an internal natural isomorphism. We call $(g, r, j, \beta)$ an \emph{algebraic trivial cofibration}. We define a morphism of algebraic trivial cofibrations to be commutative squares that preserve this structure. These form into a category $\mathbf{AlgTrivCof}$.  
\end{define}

For a $\mathbb{TC}$-coalgebra $(g, \alpha)$, define $r: \Y \to \A$ as the following composite:

\begin{equation*}
    \begin{tikzcd}
        \Y \arrow[r, "\alpha"] & \mathbb{E}(\mathbf{W}g) \arrow[r, "\mathbf{TFW}g"] & \Map(g) \arrow[r] & \A \times \Y \arrow[r] & \A.
    \end{tikzcd}
\end{equation*}

and define $\beta$ as the following composite in $\E$:

\begin{equation*}
    \begin{tikzcd}
        Y_0\arrow[r, "\alpha_0"] & E(\mathbf{W}g)_0 \arrow[r, "\mathbf{TFW}g_0"] & \Map(g)_0 \arrow[r] & Y_1.
    \end{tikzcd}
\end{equation*}

Moreover, $\alpha$ exhibits $g_0$ as a retract of a complement inclusion in $\CatE^{\mathbf{2}}$, and so $g_0$ is a complemented inclusion on objects with isomorphism $$j:= (f_0,\iota_{\Map(g)_0}^*(\alpha)): A_0 +\Map(g)_0 \times_{A_0 + \Map(g)_0} Y_0  \to Y_0.$$

In \cref{lem: coalg to retract}, we show that $(g, r, j, \beta)$ is an algebraic trivial cofibration. To prove this, it will be useful to utilise the following concept. 

\begin{define}
\label{def: algeso}
    Let $g: \A \to \Y$ be a functor with $s: Y_0 \to \Map(g)_0$  a splitting of the map 
    \begin{equation*}
        \begin{tikzcd}
            \Map(g)_0 \arrow[r] & Y_1 \arrow[r, "d_0"] & Y_0.
        \end{tikzcd}
    \end{equation*}

We call $(g,s)$ \emph{algebraically essentially surjective on objects}. 
\end{define}

\begin{lem}
\label{lem: algeso}
   Let $(g, \alpha)$ be a $\mathbb{TC}$-coalgebra. Then $(g, \mathbf{TFW}(g)_0\alpha_0)$ is algebraically essentially surjective on objects.  
\end{lem}

\begin{proof}
    This is witnessed by the following commutative diagram.

    \begin{equation*}
        \begin{tikzcd}
            Y_0 \arrow[r, "\alpha"] \arrow[ddd, equal] & E(\mathbf{W}g)_0 \arrow[dddl, "Fg_0"'] \arrow[r, "\mathbf{TFW}g_0"] & \Map(g)_0 \arrow[dl] \arrow[d] \\
            & A_0 \times Y_0 \arrow[ddl] \arrow[d, "g_0 \times 1"] & Y_1 \arrow[dd, "d_0"] \arrow[dl] \\
            & Y_0 \times Y_0 \arrow[dr] & \\
            Y_0 \arrow[rr, equal] & & Y_0
        \end{tikzcd}
    \end{equation*}
\end{proof}

\begin{lem}
\label{lem: coalg to retract}
    Let $(g, \alpha)$ be a $\mathbb{TC}$-coalgebra. Then $(g, r, j, \beta)$ as defined above is an algebraic trivial cofibration.
\end{lem}

\begin{proof}
    We first show that $r$ is a retract of $g$. This is witnessed by the following diagram.
\begin{equation*}
    \begin{tikzcd}
        \A \arrow[r, "g"] \arrow[rr, bend right, "\mathbf{TC}g"] \arrow[ddrr, "\mathbf{W}g", bend right] \arrow[dddrr, "(1 {,} g)", bend right]  \arrow[dddd, equal] & \Y \arrow[r, "\alpha"] & \mathbb{E}(\mathbf{W}g) \arrow[dd, "\mathbf{TFW}g"] \\
        & & \\
        & & \Map(g) \arrow[d] \\
        & & \A \times \Y \arrow[d] \\
        \Y \arrow[rr, equal] & & \Y
    \end{tikzcd}
\end{equation*}

Next, we show that $\beta: Y_0 \to Y_1$ is an internal natural isomorphism $gr \Rightarrow \text{id}_{\Y}.$  Firstly, the diagram

\begin{equation*}
    \begin{tikzcd}
        Y_0\arrow[r, "\alpha_0"] \arrow[rrrd, equal] & E(\mathbf{W}g)_0 \arrow[r, "\mathbf{TFW}g_0"] & \Map(g)_0 \arrow[r] & Y_1 \arrow[d, "d_0"] \\
        & & & Y_0
    \end{tikzcd}
\end{equation*}

commutes because of \cref{lem: algeso}.

The following commutes by definition of the constructions involved.

\begin{equation*}
    \begin{tikzcd}
        Y_0 \arrow[r, "\alpha_0"] \arrow[ddd, "r_0"']& E(\mathbf{W}g)_0  \arrow[r, "\mathbf{TFW}g_0"] & \Map(g)_0 \arrow[dd]  \arrow[dl] \\
        & A_0 \times Y_0 \arrow[d, "g_0 \times 1"] \arrow[ddl] &  \\
        & Y_0 \times Y_0 \arrow[dr] & Y_1 \arrow[d, "d_1"] \arrow[l] \\
       A_0 \arrow[rr, "g_0"'] & & Y_0
    \end{tikzcd}
\end{equation*}

The other conditions for natural isomorphism follow easily. 
\end{proof}

For an algebraic trivial cofibration $(g, r, j, \beta)$, we can define a coalgebra structure on $g$ by the maps which follow. 

Note that in $\CatE$, a natural isomorphism $\beta:  gr \Rightarrow \text{id}$ corresponding to a morphism $\overline{\beta}: Y_0 \to Y_1$ in $\E$ can be written as a functor $\underline{\beta} : \Y \to \Y^{\mathcal{I}}$ with $\underline{\beta}_0 = \overline{\beta}$ and $\underline{\beta}_1$ defined as below:

\begin{equation*}
    \begin{tikzcd}
        Y_1 \arrow[rrd, bend left] \arrow[ddr, bend right] \arrow[dr, dashed, "\underline{\beta}_1"] & & \\
        & \Y^{\mathcal{I}}_1 \arrow[r] \arrow[d] \arrow[dr, phantom, "\lrcorner", very near start] & \iso(Y)_1 \times_{Y_0} Y_1 \arrow[d, "m"] \\
        & Y_1 \times_{Y_0} \iso(Y)_1 \arrow[r, "m"'] & Y_1
    \end{tikzcd}
\end{equation*}

In which the two maps $Y_1 \to Y_2$ are given by $(\overline{\beta} \cdot d_1, 1_{Y_1})$ and $(d_0 \cdot g_1r_1, \overline{\beta})$.  The outside of the square commutes by the axioms for an internal natural isomorphism.

Without loss of generality, we can assume that on objects $g$ is a coproduct inclusion $A_0 \hookrightarrow A_0 + C$. First we define the following maps.

\begin{equation*}
    \begin{tikzcd}
        \Y \arrow[rrd, bend left, "\underline{\beta}"] \arrow[ddr, bend right, "(r{,}1)"']  \arrow[dr, dashed, "\hat{\beta}"] & & \\
        & \Map(g) \arrow[r] \arrow[d] \arrow[dr, phantom, "\lrcorner", very near start] & \Y^{\mathcal{I}} \arrow[d] \\
       & \A \times \Y \arrow[r] & \Y \times \Y
    \end{tikzcd}
\end{equation*}

Then we define $\tau: C \to A_0 + \Map(g)_0$ as the following composite:

\begin{equation*}
    \begin{tikzcd}
        C \arrow[r, hook] & A_0 + C \arrow[r, "\hat{\beta}_0"] & \Map(g)_0 \arrow[r, hook] & A_0 + \Map(g)_0.
    \end{tikzcd}
\end{equation*}

We define $\alpha^*: \Y \to \mathbb{E}(\mathbf{W}g)$ on objects:

\begin{equation*}
    \begin{tikzcd}
    & A_0 + \Map(g)_0 & \\
        A_0 \arrow[r, hook] \arrow[ur, "\TC g_0", hook, bend left] & A_0 + C \arrow[u, dashed, "\alpha_0^*"] & C \arrow[l, hook'] \arrow[ul, "\tau"', bend right]
    \end{tikzcd}
\end{equation*}

\begin{lem}
\label{eq: TFWgalpha = beta}
   Let $(g, r, j, \beta)$  be an algebraic trivial cofibration and consider the maps $\hat{\beta}$ and $\alpha^*_0$ as defined above. Then we have $\mathbf{TFW}(g_0)\alpha_0^* = \hat{\beta}_0$
\end{lem}

\begin{proof}
    This is witnessed by the following diagrams, using the universal properties of the coproduct and the pullback.

\begin{equation*}
    \begin{tikzcd}
        & \Map(g)_0 & \\
        & A_0 + \Map(g)_0 \arrow[u, "\mathbf{TFW}g_0"] & \Map(g)_0 \arrow[l, hook'] \arrow[ul, equal, bend right] \\
        A_0 \arrow[r, hook] \arrow[ur, bend left, hook] \arrow[uur, bend left, "\mathbf{W}g"] & A_0 + C \arrow[u, "\alpha^*"] & C \arrow[l, hook'] \arrow[ul, "\tau"] \arrow[u, "\hat{\beta} \iota_C"] 
    \end{tikzcd}
\end{equation*}

in which the left side commutes due to the commutativity of the following diagram:

\begin{equation}
\label{eq: betag}
    \begin{tikzcd}
    \A \arrow[dr, "g"] \arrow[dddrr, "(1 {,} g)"', bend right] \arrow[rr, "i"] & &  \A^{\mathcal{I}} \arrow[ddr, "g"] & \\ 
        & \Y \arrow[rrd, bend left, "\underline{\beta}"'] \arrow[ddr, bend right, "(r{,}1)"]  \arrow[dr, dashed, "\hat{\beta}"] & & \\
       &  & \Map(g) \arrow[r] \arrow[d] \arrow[dr, phantom, "\lrcorner", very near start] & \Y \arrow[d] \\
      &  & \A \times \Y \arrow[r] & \Y \times \Y
    \end{tikzcd}
\end{equation}

and the right ride commutes due to the definition of $\tau$.

\end{proof}

On morphisms, we define 

\begin{equation*}
    \begin{tikzcd}
        Y_1 \arrow[rrd, bend left, "\hat{\beta}_1"] \arrow[dd, "(d_1{,}d_0)"']  \arrow[dr, dashed, "\alpha_1^*"] & & \\
        & E(\mathbf{W}g)_1 \arrow[r] \arrow[d] \arrow[dr, phantom, "\lrcorner", very near start] & \Map(g)_1 \arrow[d] \\
       Y_0 \times Y_0 \arrow[r, "\alpha_0^* \times \alpha_0^*"'] & E(\mathbf{W}g)_0 \times E(\mathbf{W}g)_0 \arrow[r] & \Map(g)_0 \times \Map(g)_0
    \end{tikzcd}
\end{equation*}

in which the outside of the diagram commutes as witnessed by \cref{eq: TFWgalpha = beta} and the fact that $\hat{\beta}$ is an internal functor. It is not difficult to show that $\alpha^* : = (\alpha_0^*, \alpha_1^*)$ assembles into an internal functor $\Y \to \mathbb{E}(\mathbf{W}g)$. 

 \begin{lem}
 \label{lem: algret to coalg}
     Let $(g, r, j, \beta)$ be an algebraic trivial cofibration. Then $(g, \alpha^*)$ is a $\mathbb{TC}$-coalgebra. 
 \end{lem}

 \begin{proof}
    Without loss of generality, we can assume that on objects $g$ is a coproduct inclusion $g_0: A_0 \hookrightarrow A_0 + C$. We must show that the following diagram commutes: 

     \begin{equation*}
         \begin{tikzcd}
             \A \arrow[r, "\TC g"] \arrow[d, "g"'] & \mathbb{E}(\mathbf{W}g) \arrow[d, "\mathbf{F}g"] \\
             \Y \arrow[r, equal] \arrow[ur, "\alpha^*"] & \Y.
         \end{tikzcd}
     \end{equation*}

     On objects, the top triangle commutes by definition of $\alpha^*$ and the bottom diagram commutes by the following two diagrams.

     \begin{equation*}
        \begin{tikzcd}
            A_0 \arrow[r, hook] \arrow[dd, hook]  \arrow[rr, bend right, "\mathbf{TC}g_0"] \arrow[ddrr, "g"', bend right = 10] & A_0 + C \arrow[r, "\alpha_0^*"]   & A_0 + \Map(g)_0 \arrow[dd, "\mathbf{F}g_0"] \\
            & & \\
            A_0 + C \arrow[rr, equal] & & A_0 + C.
        \end{tikzcd}
    \end{equation*} 

 \begin{equation*}
        \begin{tikzcd}
            C \arrow[r, hook] \arrow[ddd, hook]   & A_0 + C \arrow[rrr, "\alpha_0^*"] \arrow[dr, "\hat{\beta}"]  \arrow[ddrr, "(r{,}1)"', bend right ]& & &  A_0 + \Map(g)_0 \arrow[ddd, "\mathbf{F}g_0"] \arrow[dl, "\mathbf{TFW}g_0"] \\
            & & \Map(g)_0 \arrow[r, equal]  \arrow[urr, hook] & \Map(g)_0 \arrow[d] & \\
            & & & A_0 \times (A_0 + C) \arrow[dr] & \\
            A_0 + C \arrow[rrrr, equal] & & & & A_0 + C.
        \end{tikzcd}
    \end{equation*} 
    
On morphisms, the following diagram witnesses that $\alpha_1^* g_1 = \mathbf{TC}g$.

\begin{equation*}
    \begin{tikzcd}[column sep = large]
        A_1 \arrow[rrr, "\mathbf{W}g_1", bend left] \arrow[r, "g_1"]  \arrow[d] & Y_1 \arrow[rr, "\hat{\beta}_1"', bend left]\arrow[d] \arrow[r, "\alpha^*_1"] & E(\mathbf{W}g)_1 \arrow[d] \arrow[r, "\mathbf{TFW}g_1"] & \Map(g)_1 \arrow[d] \\
        A_0 \times A_0 \arrow[r, "g_0 \times g_0"'] & Y_0 \times Y_0 \arrow[r, "\alpha^*_0 \times \alpha^*_0"'] & E(\mathbf{W}g)_0 \times E(\mathbf{W}g)_0 \arrow[r, "\mathbf{TFW}g_0 \times \mathbf{TFW}g_0"'{yshift=-5pt}] & \Map(g)_0 \times \Map(g)_0
    \end{tikzcd}
\end{equation*}

In this diagram, the top triangle commutes by \cref{eq: TFWgalpha = beta}. The outer square for this diagram is the defining square for $\mathbf{TC}g_1.$ Finally, the following diagram shows that $\mathbf{F}(g)_1 \alpha_1^* = 1_{Y_1}.$ 

\begin{equation*}
    \begin{tikzcd}
        Y_1 \arrow[r, "\alpha_1^*"] \arrow[d, equal] \arrow[rr, "\hat{\beta}_1", bend right] & E(\mathbf{W}g)_1 \arrow[r] & \Map(g)_1 \arrow[r] & A_1 \times Y_1 \arrow[d] \\
        Y_1 \arrow[rrr, equal] & & & Y_1 
    \end{tikzcd}
\end{equation*}

 \end{proof}

 These processes are mutually inverse, which can be checked using the universal properties involved. 

 \begin{prop}\label{prop: alg triv cof}
     There is an isomorphism of categories $\mathbb{TC}\text{-}\mathbf{Coalg} \cong \mathbf{AlgTrivCof}$.
 \end{prop}

We have therefore shown that we can translate between $\mathbb{TC}$-coalgebras and algebraic trivial cofibrations.

Next, we turn our attention to fibrations. 

\begin{define}
    Let $f: \X \to \Y$ be an internal functor in $\CatE$ and $k: X_0 \times_{Y_0} \iso(\Y)_0 \to X_1$ a splitting of the map $(d_0, f_1): \iso(\X)_0 \to X_0 \times_{Y_0} \iso(Y)_0$. We call $(f,k)$ a \emph{cloven isofibration.} 
    
    Let $(f: \X \to \Y , k)$ and $(g: \A \to \B, l)$ be cloven isofibrations. We define a \emph{morphism of cloven isofibrations} $(u,v): (f,k) \to (g,l)$ to be a morphism $(u,v): f \to g$ in $\CatE^{\mathbf{2}}$ such that the following diagram commutes.

    \begin{equation*}
        \begin{tikzcd}[column sep = large]
            X_0 \times_{Y_0} \iso(\Y)_0 \arrow[r, "u_0 \times_{v_0} \iso(v)_0"] \arrow[d, "k"'] & A_0 \times_{B_0} \iso(\B)_0 \arrow[d, "l"] \\
            \iso(\X)_0 \arrow[r, "\iso(u)_0"'] & \iso(\A)_0.
        \end{tikzcd}
    \end{equation*}

    Denote the category of cloven isofibrations and morphisms of cloven isofibrations by $\mathbf{ClovenIsofib}.$
\end{define}

Note that by \cref{rem: simplified expressions}, $X_0 \times_{Y_0} \iso(Y)_0 \cong \Map(f)_0$. 

By \cref{prop: characterisation of isofibrations}, if $(f,k)$ is a cloven isofibration, then $f$ is an isofibration, with $k$ providing a proof of this fact. We show that we can uniformly translate between $\mathbb{F}$-algebras and cloven isofibrations. 

Let $(f: \X \to \Y, \phi: \mathbb{E}(\mathbf{WF}f) \to \X)$ be an $\mathbb{F}$-algebra. Define the following maps: $$h:=(\text{id}_{X_0} \times d_1 \pi_{Y_1}, \pi_{Y_1}): X_0 \times_{Y_0}Y_1 \to \Map(f)_1$$  $$b:=  (f_1 \cdot i_{\X}\times_{Y_0}\text{id}_{Y_1}, f_1 \cdot i_{\X}\times_{Y_0}\text{id}_{Y_1}): X_0 \times_{Y_0} Y_1 \to Y_2 \times_{Y_1} Y_2 $$ $$c: = (i_{\X} \times \pi_{Y_1} , b):  X_0 \times_{Y_0} Y_1 \to \Map(f)_1$$ $$s:= ((\iota_{X_0}\pi_{X_0},h), c): X_0 \times_{Y_0} Y_1 \to \mathbb{E}(\mathbf{WF}f)_1$$
$$k : = \phi_1 s : X_0 \times_{Y_0} Y_1 \to Y_1.$$

\begin{lem}\label{lem: fib to cloven isofib}
    Let $(f: \X \to \Y, \phi: \mathbb{E}(\mathbf{WF}f) \to \X)$ be an $\mathbb{F}$-algebra. Then $(f,k)$ is a cloven isofibration. 
\end{lem}

Conversely, starting with $(f: \X \to \Y,k)$ a cloven isofibration, we can can form the structure of an $\mathbb{F}$-algebra via the following maps.

Firstly, define $\phi_0: (\text{id}_{X_0}, d_1 \cdot k): E(\mathbf{W}f)_0 : = X_0 + \Map(f)_0 \to X_0.$ Next, we note that by \cref{rem: simplified expressions} in order to define a map $\phi_1: E(\mathbf{W}f)_1 \to X_1,$ it is enough to define maps $a: (X_0 \times X_0)\times_{\Map(f)_0 \times \Map(f)_0} \Map(f)_1 \to X_1$, $b: X_0 \times_{\Map(f)_0}\Map(f)_1 \to X_1$ and $c: \Map(f)_1\to X_1$. We define these maps as follows:

\begin{equation*}
    a : = \left(\begin{tikzcd}
          (X_0 \times X_0)\times_{\Map(f)_0\times \Map(f)_0} \Map_1 \arrow[r, "\pi_{\Map(f)_1}"] & \Map(f)_1 \arrow[r, "\pi_{X_1}"] & X_1
    \end{tikzcd}\right),
\end{equation*}
\begin{equation*}
  \begin{tikzcd}
    \Map(f)_1 \arrow[dr, dashed, "x"] \arrow[d, "d_0"'] \arrow[drr, "\pi_{X_1}", bend left] &  &  \\
    \Map(f)_0 \arrow[d, "k"] &   X_2 \arrow[r] \arrow[d] \arrow[dr, phantom, "\lrcorner", very near start] & X_1 \arrow[d, "d_0"] \\
     \iso(\X)_1 \arrow[r, hookrightarrow] & X_1 \arrow[r, "d_1"] & X_0 
    \end{tikzcd}
\end{equation*}
\begin{equation*}
    b : = \left(\begin{tikzcd}
         X_0 \times_{\Map(f)_0} \Map(f)_1 \arrow[r] & \Map(f)_1 \arrow[r, "x"] & X_2 \arrow[r, "m"] & X_1 
    \end{tikzcd} \right),
\end{equation*}
\begin{equation*}
    c := \left(\begin{tikzcd}[column sep = huge]
         \Map(f)_1 \arrow[r, "(\pi_{X_1}{,} d_1{,} d_0)"] & X_1 \times \Map(f)_0 \times \Map(f)_0 \arrow[r, "(\text{id}{,} \text{inv} \cdot k{,} k)"] & X_3 \arrow[r, "m^2"] & X_1
    \end{tikzcd} \right).
\end{equation*}

We define $\phi_1: = a + b + c: E(\mathbf{W}f)_1 \to X_1.$ It is not hard to see that we have an internal functor $\phi : = (\phi_0, \phi_1) : \mathbb{E}(\mathbf{W}f) \to \X.$

\begin{lem}
    Let $(f: \X \to \Y,k)$ a cloven isofibration. Then $(f, \phi)$ is an $\mathbb{F}$-algebra. 
\end{lem}

\begin{proof}
    It is clear that on objects $\phi_0 \cdot \mathbf{TC}f_0 = \text{id}_{X_0}$ and the following commutative diagram completes the proof that $f_0 \cdot \phi_0 = \mathbf{F}f_0.$

    \begin{equation*}
        \begin{tikzcd}
            \Map(f)_0 \arrow[r, "k"] \arrow[dr, equal]  & \iso(\X)_0 \arrow[r, "d_1"] \arrow[d, "(d_1{,} \iso(f)_0)"] & X_0 \arrow[dd, "f_0"] \\
             & \Map(f)_0 \arrow[d] & \\
             & X_0 \times Y_0 \arrow[r, "\pi_{Y_0}"'] & Y_0.
        \end{tikzcd}
    \end{equation*}
    We note that $\mathbf{TC}f_1$ factors through $(X_0 \times X_0)\times_{\Map(f)_0\times \Map(f)_0} \Map_1$ and $a\cdot\mathbf{TC}f = 1_{\X_1}$ by construction. It remains to show that $f_1 \cdot b = \pi_{Y_1} \pi_{\Map(f)_1}$ and $f_1 \cdot c = \pi_{Y_1}$ as these morphisms are $\mathbf{F}f_1|_{X_0 \times_{\Map(f)_0} \Map(f)_1}$ and $\mathbf{F}f_1|_{\Map(f)_1}$ respectively. This can be shown representably.
\end{proof}

These processes are mutually inverse, which can again be checked using the universal properties involved. 

\begin{prop}\label{prop: cloven isofib}
    There is an isomorphism of categories $\mathbb{F}\text{-}\mathbf{Alg} \cong \mathbf{ClovenIsofib}$.
\end{prop}

We also have equivalent characterisations of $\mathbb{C}$-coalgebras and $\mathbb{TF}$-algebras; for these, we omit the proof since it uses the same strategy as previously used and their construction will not be needed in the subsequent work. We record these below.

\begin{define}
    Let $f: \X \to \Y$ and an isomorphism $j: X_0 + C \to Y_0$ such that $j \cdot \iota_{X_0} = f_0$. We call $(f, j)$ an \emph{algebraic complemented inclusion on objects.}

    Let $(f: \X \to \Y,j: Y_0 \to X_0 + C), (g: \A \to \B, t: Y_0 \to D)$ be a pair of algebraic complemented inclusion on objects. A morphism $(u,v); f \to g$ is called a morphism of algebraic complemented inclusion on objects if the following diagram commutes.

    \begin{equation*}
        \begin{tikzcd}
            X_0 \arrow[d, hookrightarrow] \arrow[r, "u_0"] & A_0 \arrow[d, hookrightarrow] \\
            X_0 + C \arrow[r, "t\cdot v_0 \cdot j^{-1}"'] & A_0 + D.
        \end{tikzcd}
    \end{equation*}

We denote the category of algebraic complemented inclusion on objects and the morphisms of these by $\mathbf{AlgCompIncObj}.$
\end{define}

\begin{prop}
    There is an isomorphism of categories $\mathbb{C}\text{-}\mathbf{Coalg}\cong\mathbf{AlgCompIncObj}$.
\end{prop}

\begin{define}
    Let $f: \X \to \Y$,  $s: \Y \to \X$ be a splitting of $f$ and $\beta: \text{id}_{\X} \Rightarrow sf$ a natural isomorphism. We call $(f,s, \beta)$ an \emph{algebraic split epi equivalence}. We define a morphism of algebraic split epi equivalences to be commutative squares that preserve the choice of splitting. These form into a category $\mathbf{AlgSplitEpiEq}$.  
\end{define}

\begin{prop}
    There is an isomorphism of categories $\mathbb{TF}\text{-}\mathbf{Alg} \cong \mathbf{AlgSplitEpiEq}$.
\end{prop}

\section{Type theoretic aspects}\label{sec: TTAWFS}

In this section, we restrict our attention to the setting of $\E$ a locally cartesian closed locos with coequalisers, which is a locally cartesian closed and lextensive category with parametrised list objects and coequalisers. Such an $\E$ gives us an example of a lextensive cartesian closed category pullback stable coequalisers in which the forgetful functor $\mathcal{U}: \CatE \to \mathbf{Gph}(\E)$ has left adjoint \cite{maietti2010joyal, Hughes2024Colimits}. We also restrict our attention from $\CatE$ to $\GpdE$--- this is because one ingredient we need to model type theory is exponentiability of isofibrations. It is not true that all isofibrations in $\CatE$ are exponentiable; however, this is true for isofibrations in $\GpdE$. In fact, in \cref{lem: lcc} we show that that isofibrations in $\GpdE$ are exponentiable if and only if $\E$ is locally cartesian closed, which shows that we must assume local cartesian closure. For $\E$ a locally cartesian closed locos with coequalisers, $\GpdE$ has the structure of an algebraic model structure by \cref{cor: algebraic model structure on GpdE}.

The aim of this section is to prove that the $(\mathbb{TC}, \mathbb{F})$ algebraic weak factorisation system on $\GpdE$ has the structure of a type theoretic algebraic weak factorisation system in the sense of \cite[Definition 4.10]{gambino2023models} which we recall in \cref{def:ttawfs}. By \cite[Theorem 3.12]{gambino2023models}, this proves that the cloven isofibrations in $\GpdE$ give us an algebraic model of MLTT with strictly stable choices of of $\Sigma\text{-}, \Pi\text{-}$  and $\texttt{Id}$-types.

\subsection{Exponentiability}

Recall that a morphism $f: X \to Y$ in a category $\C$ is called exponentiable if the pullback functor $f^*: \C/Y \to \C/X$ has right adjoint. We say that an algebraic weak factorisation system $(\mathbb{L}, \mathbb{R}) $ satisfies the \emph{exponentiability property} if any map in the image of the forgetful functor $\mathbb{R}\text{-}\mathbf{Alg} \to \C^{\mathbf{2}}$ is exponentiable \cite[Definition 3.4]{gambino2023models}. Exponentiability of isofibrations for internal groupoids is studied by Niefield-Pronk in \cite{niefield2019internal}, and immediately gives us the result of \cref{prop: exp}. 
\begin{prop}
\label{prop: exp}
    $(\mathbb{TC}, \mathbb{F})$ satisfies the exponentiability property.
\end{prop}

\begin{proof}
    By \cref{prop: algebras for pointed endofunctors}, $\mathbb{F}$-algebras are precisely the isofibrations with a choice of lifting. Given that $\E$ is locally cartesian closed, we can apply Theorem 4.5 of \cite{niefield2019internal}. 
\end{proof}

\subsection{Frobenius structure}
\label{sec: Frob}

The aim of this section is to construct a functorial Frobenius structure on $(\mathbb{TC}, \mathbb{F})$, as recalled in the following definition. 

\begin{define}[\cite{gambino2023models}, Definition 3.5]
    Let $(\mathbb{L}, \mathbb{R})$ be an algebraic weak factorisation system on a category $\C$. A \emph{functorial Frobenius structure} on $(\mathbb{L}, \mathbb{R})$ is given by the dotted map in the following commutative diagram:

    \begin{equation*}
       \begin{tikzcd}
            \mathbb{R}\text{-}\mathbf{Alg} \times_{\C} \mathbb{L}\text{-}\mathbf{Coalg} \arrow[d] \arrow[r, dashed, "\widehat{\mathbf{pb}}"] & \mathbb{L}\text{-}\mathbf{CoAlg} \arrow[d] \\
            \C^{\mathbf{2}} \times_{\C} \C^{\mathbf{2}} \arrow[r, "\mathbf{pb}"'] & \C^{\mathbf{2}}.
       \end{tikzcd}
    \end{equation*}

    in which $\mathbf{pb}: \C^{\mathbf{2}} \times_{\C} \C^{\mathbf{2}} \to \C^{\mathbf{2}}$ denotes pullback: $(f,g) \mapsto f^*(g)$. 
\end{define}

To provide a functorial Frobenius structure in our setting, we must construct a functor $\widehat{\mathbf{pb}}: \mathbb{F}\text{-}\mathbf{Alg} \times_{\CatE} \mathbb{TC}\text{-}\mathbf{Coalg} \to \mathbb{TC}\text{-}\mathbf{Coalg}$. To do this, we show that algebraic trivial cofibrations are closed under pullback along cloven isofibrations in a uniform way. By \cref{prop: cloven isofib,prop: alg triv cof}, the result follows.

First, note that in an extensive category, complemented inclusions are stable under pullback along any map by definition \cite{carboni1993introduction}. Given $i: A \to A+C$ and $f: X \to A+C$ we can deduce that $X \cong A\times_{A+C}X + C\times_{A+C}X$. Therefore, given an algebraic trivial cofibration $(g, r, j, \beta)$ and a cloven isofibration $(f, k)$, $f^*(g)$ has the structure of a complemented inclusion on objects functor, which we denote by $j^*$. Next, we show that retract equivalences are closed under pullback along cloven isofibrations in a canonical way.

Let $(g: \A \to \Y, r, \beta)$ be a retract equivalence and $(f: \X \to \Y, k: X_0 \times_{Y_0} Y_1 \to X_1)$ be a cloven isofibration. We define a map $r^*_0: X_0 \to A_0 \times_{Y_0} X_0$ by the following:

\begin{equation*}
    \begin{tikzcd}
        X_0 \arrow[ddr, bend right, equal] \arrow[r, "f_0"] \arrow[dr, dashed, "\beta'"] & Y_0 \arrow[dr, "\overline{\beta}"] & \\
        & Y_1 \times_{Y_0} X_0 \arrow[r] \arrow[d] \arrow[dr, phantom, "\lrcorner", very near start] & Y_1 \arrow[d, "d_1"] \\
        & X_0 \arrow[r, "f_0"'] & Y_0 
    \end{tikzcd}
\end{equation*}
\begin{equation*}
    \begin{tikzcd}
        X_0 \arrow[ddr, bend right, "d_0k\beta'"'] \arrow[r] \arrow[dr, dashed, "r^*_0"] & Y_0 \arrow[dr, "r_0"] & \\
        & X_0 \times_{Y_0} A_0 \arrow[r] \arrow[d] \arrow[dr, phantom, "\lrcorner", very near start] & A_0 \arrow[d, "g_0"] \\
        & X_0 \arrow[r, "f_0"'] & Y_0 
    \end{tikzcd}
\end{equation*}

where the outside of the second diagram commutes by the following diagram.

\begin{equation*}
    \begin{tikzcd}
        X_0 \arrow[r, "\beta'"] \arrow[dd, "f_0"'] & X_0\times_{Y_0} Y_1 \arrow[r, "k"] \arrow[d, equal] & X_1 \arrow[r, "d_0"] \arrow[dl, "l"] \arrow[d, "f_1"] & X_0 \arrow[dd, "f_0"] \\
        & X_0 \times_{Y_0}Y_1 \arrow[r] &  Y_1 \arrow[dr, "d_0"] & \\
        Y_0 \arrow[r, "r_0"']\arrow[rru, "\overline{\beta}"] & A_0 \arrow[rr, "g_0"'] & & Y_0
    \end{tikzcd}
\end{equation*}

Next, recall that the model structure on $\GpdE$ is right proper since all objects are fibrant, and therefore weak equivalences are closed under pullback along fibrations \cite{hirschhorn2003model}. Therefore $f^*(g)$ is a weak equivalence and is in particular fully faithful. We construct a morphism $r^*_1: X_1 \to A_1 \times_{Y_1} X_1$ as follows:

\begin{equation*}
    \begin{tikzcd}[column sep = large]
        X_1 \arrow[dr, dashed, "r^*_1"] \arrow[drr, bend left, equal] \arrow[dd] & & \\
        &  A_1 \times_{Y_1} X_1 \arrow[d] \arrow[r, "f^*(g)_1"] \arrow[dr, "\lrcorner", phantom, very near start] & X_1 \arrow[d] \\
        X_0 \times X_0 \arrow[r, "r^*_0 \times r^*_0"'] & A_0 \times_{Y_0} X_0 \times A_0 \times_{Y_0} X_0 \arrow[r, "f^*(g)_0 \times f^*(g)_0"'] & X_0 \times X_0.
    \end{tikzcd}
\end{equation*}

We define an internal functor by $r^* = (r^*_0, r^*_1).$

Define a map $\overline{\beta^*}: X_0 \to X_1$ as the following composite:

\begin{equation*}
    \begin{tikzcd}
        X_0 \arrow[r, "\beta'"] & X_0 \times_{Y_0} Y_1 \arrow[r, "k"] & X_1.
    \end{tikzcd}
\end{equation*}

Consider the following diagrams:

\begin{equation*}
    \begin{tikzcd}[column sep = huge]
        X_0 \arrow[r, "\beta'"] \arrow[ddrr, equal, bend right] \arrow[dr, "\beta'"']& X_0 \times_{Y_0} Y_1 \arrow[r, "k"] \arrow[d, equal] & X_1  \arrow[dl, "l"] \arrow[dd, "d_1"] \\
        & X_0 \times_{Y_0} Y_1 \arrow[dr] & \\
        & &  X_0
    \end{tikzcd}
\end{equation*}

\begin{equation*}
    \begin{tikzcd}[column sep = huge]
        X_0 \arrow[r, "\beta'"] \arrow[d,  "r^*_0"'] & X_0 \times_{Y_0} Y_1 \arrow[r, "k"]  & X_1   \arrow[d, "d_0"] \\
      X_0 \times_{Y_0} A_0 \arrow[rr, "f^*(g)_0"'] & &  X_0
    \end{tikzcd}
\end{equation*}

in which the latter commutes by definition of $\beta'.$

It is not hard to prove that $\overline{\beta^*}$ defines a natural isomorphism $\beta*: f^*(g) r^* \Rightarrow 1$. We conclude that algebraic trivial cofibrations are closed under pullback along isofibrations,  the details of which are contained in the construction.

\begin{lem}
\label{lem: algret closed under pullback along fibs}
    Let $(g: \A \to \Y, r, j, \beta)$ be an algebraic trivial cofibration and $(f,k)$ be a cloven isofibration. Then $(f^*(g), r^*, j^*, \beta^*)$ is an algebraic trivial cofibration. 
\end{lem}

We are now able to prove the following result.

\begin{prop}
\label{prop: Frob}
      There is a Frobenius structure on $(\mathbb{TC}, \mathbb{F})$. 
\end{prop}

\begin{proof}
     Let $(g: \A \to \Y, \alpha)$ be a $\mathbb{TC}$-coalgebra and $(f: \X \to \Y,s)$ be an $\mathbb{F}$-algebra. We uniformly construct a $\mathbb{TC}$-coalgebra structure on the pullback $f^*(g).$ By \cref{lem: coalg to retract}, we translate $(g, \alpha)$ into a algebraic trivial cofibration $(g, r, j, \beta)$ and by \cref{lem: fib to cloven isofib}, we translate $(f, s)$ into a cloven isofibration $(f,k)$. Therefore, by \cref{lem: algret closed under pullback along fibs} $(f^*(g), r^*, j^*, \beta^*)$ is an algebraic trivial cofibration. Whence, we can apply \cref{lem: algret to coalg} to obtain a $\mathbb{TC}$-coalgebra $(f^*(g), \alpha^*)$, as required. 
\end{proof}

\subsection{Stable functorial choice of path objects}
\label{sec: sfcpo}

The aim of this section is to show that $(\mathbb{TC}, \mathbb{F})$ has a stable functorial choice of path objects. For $f: X \to Y$ in $\C$, denote the diagonal map by $\Delta_f: X \to X\times_Y X$. This extends to a functor $\Delta_{-}: \C^{\mathbf{2}} \to \C^{\mathbf{2}}.$ We recall the following definitions below from \cite{gambino2023models}.

\begin{define}
    Let $(\mathbb{L}, \mathbb{R})$ be an algebraic weak factorisation system on $\C$. A \emph{functorial factorisation of the diagonal} is a functor $\mathcal{P}: \C^{\mathbf{2}} \to \C^{\mathbf{2}}\times_{\C}\C^{\mathbf{2}}$ such that $$f: X \to Y \mapsto \begin{tikzcd}
        X \arrow[r, "\lambda_f"] & \mathcal{P}X \arrow[r, "\rho_f"] & X \times_Y X 
    \end{tikzcd}$$

    in which $\rho_f \cdot \lambda_f = \Delta_f.$

    Such a functorial factorisation of the diagonal is called \emph{stable} if $(h,k): (f \to f') \in \C^{\mathbf{2}}$ is a pullback if and only if $\rho_{(h{,}k)}: \rho_f \to \rho_{f'}$ is a pullback. 
    
    A \emph{stable functorial choice of path objects} consists of a lift of a stable functorial factorisation of the diagonal map, as shown in the following diagram.

    \begin{equation*}
        \begin{tikzcd}
           \mathbb{R}\text{-}\mathbf{Alg} \arrow[r, "\widehat{\mathcal{P}}", dashed] \arrow[d] &  \mathbb{L}\text{-}\mathbf{Coalg}\times_{\C} \mathbb{R}\text{-}\mathbf{Alg} \arrow[d] \\
           \C^{\mathbf{2}} \arrow[r, "\mathcal{P}"'] & \C^{\mathbf{2}} \times_{\C} \C^{\mathbf{2}}.
        \end{tikzcd}
    \end{equation*}
\end{define}

In $\CatE$, we claim that the functorial factorisation of the diagonal $(\mathbf{TC}(\Delta_{-}), \mathbf{F}(\Delta_{-})): \C^{\mathbf{2}} \to \C^{\mathbf{2}}\times_{\C}\C^{\mathbf{2}}$ gives us a stable functorial choice of path objects. 

\begin{prop}
    The functorial factorisation of the diagonal $(\mathbb{TC}(\Delta_{-}), \mathbb{F}(\Delta_{-}))$ is stable. 
\end{prop}

\begin{proof}
    Consider a cartesian square 

    \begin{equation}
    \label{eq: cartesian square}
        \begin{tikzcd}
            \A \arrow[r, "u"] \arrow[d, "g"'] \arrow[dr, "\lrcorner", very near start, phantom] & \X \arrow[d, "f"] \\
            \B \arrow[r, "v"'] & \Y.
        \end{tikzcd}
    \end{equation}

    Firstly, we argue that the following square is a pullback.
    \begin{equation}
    \label{eq: Map}
        \begin{tikzcd}
            \Map(g) \arrow[r, "\Map(u{,}v)"] \arrow[d] \arrow[dr, "\lrcorner", very near start, phantom] & \Map(f) \arrow[d] \\
            \B \arrow[r, "v"'] & \Y
        \end{tikzcd}
    \end{equation}

It is enough to show that this is true in $\s$ as $\Map(g)$ is a representable construction as it is built out of pullbacks, products and exponentiation, which are all representable notions. Therefore, this follows from (Proposition 4.4 of \cite{gambino2023models}).

By the pullback lemma, given that \cref{eq: Map} and \cref{eq: cartesian square} are pullbacks, it follows that the following square is a pullback:

\begin{equation}
\label{eq: map cartesian}
\begin{tikzcd}
       \A \arrow[r, "u"] \arrow[d, "\mathbf{C} g"'] \arrow[dr, phantom, very near start, "\lrcorner"] & \X \arrow[d, "\mathbf{C} f"] \\
    \Map(g) \arrow[r, "\Map(u{,}v)"] & \Map(f).
\end{tikzcd}
\end{equation}

     Next, we show that given a pullback square as in \cref{eq: cartesian square}, the following square is a pullback:

    \begin{equation}
    \label{eq: cartesian E(g)}
        \begin{tikzcd}
            \mathbb{E}(g) \arrow[r, "\mathbb{E}(u{,}v)"] \arrow[d, "\TF g"'] & \mathbb{E}(f) \arrow[d, "\TF f"] \\
            \B \arrow[r, "v"'] & \Y
        \end{tikzcd}
    \end{equation}

On objects, this is the square

\begin{equation*}
    \begin{tikzcd}
        A_0 + B_0 \arrow[d, "\TF g"'] \arrow[r, "E(u{,}v)_0"] & X_0 + Y_0 \arrow[d, "\TF f"] \\
        B_0 \arrow[r, "v"'] & Y_0.
    \end{tikzcd}
\end{equation*}

which is a pullback due to distributivity of pullbacks over coproducts in $\E$ as seen in the following calculation:

\begin{align*}
    (X_0 + Y_0) \times_{Y_0} B_0 & \cong (X_0 \times_{Y_0} B_0)+ (Y_0 \times_{Y_0} B_0) \\
    & \cong A_0 + B_0.
\end{align*}

On morphisms, consider the following diagram in $\E$:

\begin{equation}
\label{eq: pullback on morphisms}
    \begin{tikzcd}
       E \arrow[drr, bend left, "p"] \arrow[ddr, bend right, "q"'] \arrow[dr, dashed, "t"] & & \\
       & E(g)_1 \arrow[r, "E(u{,}v)_1"] \arrow[d, "\TF g_1"'] & E(f)_1 \arrow[d, "\TF f_1"] \\
       &  B_1 \arrow[r, "v_1"'] & Y_1.
    \end{tikzcd}
\end{equation}

We describe below how the dotted arrow $E \to E(g)_1$, making the above square into a pullback. First, define $k:E \to A_0 \times A_0$ by the universal property of $A_0 \times A_0$ as the pullback over $v_0 \times v_0: B_0 \times B_0 \to Y_0 \times Y_0$ and $f_0 \times f_0: X_0 \times X_0 \to Y_0 \times Y_0$, given that the square from \cref{eq: cartesian square} is a pullback:

\begin{equation*}
    \begin{tikzcd}
       E \arrow[d, "q"'] \arrow[r, "p"] \arrow[dr, dashed, "k"] & E(f)_1 \arrow[dr, "\TF f"] & \\
       B_1 \arrow[dr, "(d_1{,} d_0)"']  & A_0 \times A_0 \arrow[dr, phantom, very near start, "\lrcorner"] \arrow[r, "u_0 \times u_0"] \arrow[d, "g_0 \times g_0"'] & X_0 \times X_0 \arrow[d, "f_0 \times f_0"] \\
       &  B_0 \times B_0 \arrow[r, "v_0 \times v_0"'] & Y_0 \times Y_0.
    \end{tikzcd}
\end{equation*}

This in turn induces an arrow $t: E \to E(g)_1$ by the universal property of $E(g)_1$ as a pullback, given the commutativity of the outer square of the diagram below, which is by construction.

\begin{equation*}
    \begin{tikzcd}
       E \arrow[drr, bend left, "q"] \arrow[ddr, bend right, "k"'] \arrow[dr, dashed, "t"] & & \\
       & E(g)_1 \arrow[r, "\TF g_1"] \arrow[d] \arrow[dr, phantom, very near start, "\lrcorner"] & B_1 \arrow[d, "(d_1{,} d_0)"] \\
       & A_0 \times A_0 \arrow[r, "g_0 \times g_0"'] & B_0 \times B_0.
    \end{tikzcd}
\end{equation*}

It remains to show that $E(u,v)_1 t = p,$ which will commutativity of \cref{eq: pullback on morphisms}, therefore showing that the inner square is a pullback. We show that $E(u,v)_1 t = p$ by using the universal property of $E(f)_1$ as a pullback: we show that the maps agree on both of its projections. This is shown in the diagrams below.

\begin{equation*}
\begin{tikzcd}
     E \arrow[rr, "t"] \arrow[dd, "p"'] \arrow[dr, "q"] & & E(g)_1 \arrow[dl, "\TF g_1"] \arrow[d, "E(u{,}v)"] \\
    & B_1 \arrow[dr, "v_1"'] & E(f)_1 \arrow[d, "\TF f_1"] \\
    E(f)_1 \arrow[rr, "\TF f_1"'] & & Y_1 
\end{tikzcd}
\begin{tikzcd}
     E \arrow[rr, "t"] \arrow[dd, "p"'] \arrow[dr, "k"] & & E(g)_1 \arrow[dl, "(d_1{,} d_0)"] \arrow[d, "E(u{,}v)"] \\
    & A_0 \times A_0 \arrow[dr, "u_0 \times u_0"'] & E(f)_1 \arrow[d, "(d_1{,} d_0)"] \\
    E(f)_1 \arrow[rr, "(d_1{,} d_0)"'] & & X_0 \times X_0
\end{tikzcd}
\end{equation*}

Since pullbacks in $\CatE$ are calculated pointwise, it follows that \cref{eq: cartesian E(g)} is a pullback square.

Hence, factorising the cartesian square in \cref{eq: map cartesian}, we obtain the cartesian square

\begin{equation}
\label{eq: TFW}
        \begin{tikzcd}
            \mathbb{E}(\Weq g) \arrow[r, "\mathbb{E}(u{,}v)"] \arrow[d, "\TF \Weq g"'] \arrow[dr, phantom, "\lrcorner", very near start] & \mathbb{E}(\Weq f) \arrow[d, "\TF \Weq  f"] \\
            \Map(g) \arrow[r, "v"'] & \Map(f)
        \end{tikzcd}
    \end{equation}

By the pullback lemma, by pasting \cref{eq: TFW} and \cref{eq: Map}, we obtain the desired pullback square.

\end{proof}

\begin{prop}
     The factorisation of the awfs $(\mathbb{TC}, \mathbb{F})$ gives a stable functorial choice of path objects.
\end{prop}
   
\begin{proof}
    It remains to show that we can lift the stable functorial factorisation $(\mathbf{TC}(\Delta_{-}), \mathbf{F}(\Delta_{-})): \C^{\mathbf{2}} \to \C^{\mathbf{2}}\times_{\C}\C^{\mathbf{2}}$ to algebas, but this is obvious from construction; given an $\mathbb{F}$-algebra $(f, \alpha)$, $(\mathbb{TC}(\delta_f), \sigma_{\Delta_f}) \in \mathbf{TC}\text{-}\mathbf{Coalg}$ and $(\mathbf{F}(\Delta_f), \kappa_{\Delta_f}) \in \mathbb{F}\text{-}\mathbf{Alg}.$ 
\end{proof}

\subsection{A type theoretic algebraic weak factorisation system}

We are now able to state our main theorem. Recall the following definition.

\begin{define}[\cite{gambino2023models}, Definition 3.10]\label{def:ttawfs}
    Let $\C$ be a category. A \emph{type-theoretic algebraic weak factorisation system} consists of the following data:

    \begin{enumerate}
        \item An algebraic weak factorisation system $(\mathbb{L}, \mathbb{R})$  on $\C$ satisfying the exponentiability condition.
        \item a functorial Frobenius structure on $(\mathbb{L}, \mathbb{R}).$
        \item a stable, functorial choice of path objects on $(\mathbb{L}, \mathbb{R}).$
    \end{enumerate}
\end{define}

By \cite[Theorem 3.12]{gambino2023models}, any type theoretic algebraic weak factorisation system forms a model of MLTT with strictly stable choices of $\Sigma$-, $\Pi$- and $\texttt{Id}$-types. We have therefore proven the following.

\begin{thm}
\label{thm: ttawfs}
    Let $\E$ be a locally cartesian closed locos with coequalisers. Then the awfs $(\mathbb{TC}, \mathbb{F})$ on the category $\GpdE$ is equipped with the structure of a type theoretic awfs. 

    Consequently, the right adjoint splitting of the comprehension category associated to $(\mathbb{TC}, \mathbb{F})$ is a model of MLTT with strictly stable choices of $\Sigma$-, $\Pi$- and $\texttt{Id}$-types.   
\end{thm}

\section{$2$-categorical aspects}

In this section, we exploit the work in \cite{Hughes2024Colimits} and \cite{bourke2010codescent} to obtain an elementary model of MLTT in purely $(2,1)$-categorical terms. We provide a list of necessary and sufficient axioms on a $2$-category $\K$ that implies $\K \simeq \CatE$ for $\E$ a locally cartesian closed locos with coequalisers, a result we believe is of independent interest. Then, we show that we can restrict this to $(2,1)$-categories to get a purely $(2,1)$-categorical model of MLTT. 

We start by characterising a $2$-dimensional version of local cartesian closure. Note that $\E$ being locally cartesian closed does not imply that $\CatE$ is locally cartesian closed; this fails even in the case that $\E = \s$. Instead of asking that all morphisms in $\K$ are exponentiable, we ask only that the groupoidal isofibrations are exponentiable. We also remark that an equivalent characterisation is that we ask for discrete opfibrations to be exponentiable.

   \begin{prop}\label{lem: lcc}
       Let $\E$ have finite limits and finite colimits. Then $\E$ is locally cartesian closed if and only isofibrations in $\Gpd(\E)$ are exponentiable. 
   \end{prop}

   \begin{proof}
       Suppose $\E$ is a locally cartesian closed category with finite limits and coproducts. Then isofibrations in $\Gpd(\E)$ are exponentiable by \cite[Theorem 4.5]{niefield2019internal}. 
       
       Conversely, suppose isofibrations in $\Gpd(\E)$ are exponentiable. Then for any internal functor $f: X \to Y$ in $\E$, we note that $\disc(f): \disc(X) \to \disc(Y)$ is an isofibration in $\Gpd(\E)$; this is easy to see representably. Hence, we have a right adjoint to the functor $\disc(f)^*: \Gpd(\E)/ \disc(Y) \to \Gpd(\E)/ \disc(X)$ which we denote  by $\Pi_{\disc(f)}: \GpdE/\disc(Y) \to \GpdE/\disc(X)$. Construct the right adjoint to $f^*: \E/Y \to \E/X$ by composing this series of right adjoints $(-)_0 \cdot \Pi_{\disc(f)} \cdot \disc: \E/X \to \E/Y$ which is left adjoint to $\Pi_0 \cdot \disc(f)^* \cdot \disc: \E/Y \to \E/X$ which we can easily verify is equal to $f^*: \E/Y \to \E/X$. For $\Pi_0$ to exist, we needed to assume that we had coequalisers.
   \end{proof}

   We can state the property of isofibrations in $\GpdE$ representably for any $2$-category $\K$.

   \begin{define}\label{def: groupoidal isofibrations}
       Let $\K$ be a $2$-category. We say a morphism $F : X \to Y$ in $\K$ is a groupoidal isofibration if for all $A \in \K$ the functor $\K(A, F): \K(A, X) \to \K(A, Y)$ is an isofibration between groupoids. 
   \end{define}

   For $\K = \CatE$, this notion recovers exactly isofibrations in $\GpdE.$ 
   
   We are now able to give a complete characterisation of  locally cartesian closed locoi in purely $2$-categorical terms.
   
     \begin{thm}\label{thm: final characterisation}
       Let $\E$ be a locally cartesian closed locos with coequalisers. Then the $2$-category $\K : = \CatE$ satisfies the conditions listed below. Conversely, if $\mathcal{K}$ satisfies the conditions listed below, then there is a $2$-equivalence $\K \simeq \mathbf{Cat}\left(\E\right)$ where $\E:= \mathbf{Disc}\left(\mathcal{K}\right)$, in which $\E$ is a locally cartesian closed locos with coequalisers.

    \begin{enumerate}
        \item $\mathcal{K}$ has pullbacks and powers by $\mathbf{2}$.
        \item $\mathcal{K}$ has codescent objects of categories internal to $\mathcal{K}$ whose source and target maps form a two-sided discrete fibration.
        \item Codescent morphisms are effective in $\mathcal{K}$.
        \item Discrete objects in $\mathcal{K}$ are projective, in the sense of Definition 4.13 of \cite{bourke2010codescent}.
        \item For every object $A \in \mathcal{K}$, there is a projective object $P \in \mathcal{K}$ and a codescent morphism $c: P \rightarrow A$.
        \item $\K$ is lextensive.
        \item Groupoidal isofibrations in $\K$ are exponentiable.
        \item $\K$ has finite $2$-colimits. 
        
    \end{enumerate}
   \end{thm}

   \begin{proof}
       Suppose $\E$ is a locally cartesian closed locos with coequalisers. Then $\CatE$ satisfies (1)-(5) by \cite{bourke2010codescent}, (6) by \cite{Hughes2024Elementary},  (7) by \cite[Theorem 4.5]{niefield2019internal} and (8) by \cite[Theorem 5.2]{Hughes2024Colimits}.
       
       Conversely, suppose $\K$ is a $2$-category satisfying (1)-(8). Again, (1)-(5) shows that $\K \simeq \CatE$, (6) shows that $\E$ is lextensive by \cite[Proposition 3.3, Lemma 5.2]{Hughes2024Elementary}. (7) implies that $\E$ is locally cartesian closed by Proposition~\ref{lem: lcc}. (8) implies that $\K$ has a natural numbers object by \cite[Corollary 6.3]{Hughes2024Colimits}. Moreover, it has coeqalisers by \cite[Lemma 6.1]{Hughes2024Colimits}. In a local cartesian closed category, having a natural numbers object is equivalent to having an parametrised list objects by \cite[Theorem 2.5.17]{johnstone2002sketches}; hence $\E$ is a locally cartesian closed locos with coequalisers. 
   \end{proof}

Below, we state a $(2,1)$ version of \cite[Theorem 4.19]{bourke2010codescent}.

\begin{thm}\label{thm: bourke}
    Let $\E$ be a category with pullbacks. Then the $(2,1)$-category $\GpdE$ satisfies the conditions listed below. Conversely, if $\K$ is a $(2,1)$-category satisfying the conditions listed below, then there is a $(2,1)$-equivalence $\K \simeq \GpdE$ where $\E:= \mathbf{Disc}(\K)$. 

    \begin{enumerate}
        \item $\mathcal{K}$ has pullbacks and powers by $\mathcal{I}$.
        \item $\mathcal{K}$ has codescent objects of groupoids internal to $\mathcal{K}$ whose source and target maps form a two-sided discrete fibration.
        \item Codescent morphisms are effective in $\mathcal{K}$.
        \item Discrete objects in $\mathcal{K}$ are projective, in the sense of Definition 4.13 of \cite{bourke2010codescent}.
        \item For every object $A \in \mathcal{K}$, there is a projective object $P \in \mathcal{K}$ and a codescent morphism $c: P \rightarrow A$.
    \end{enumerate}
\end{thm}

\begin{proof}

    It is not hard to see that for any $\E$ with pullbacks, the $(2,1)$-category $\GpdE$ satisfies the required conditions. 

    Conversely, suppose that $\K$ is a $(2,1)$-category satisfying the conditions of the Theorem. Note that conditions (1)-(5) are the conditions for \cite[Theorem 4.19]{bourke2010codescent} with the exceptions of (1) and (2), which we show are $(2,1)$-versions. For (1), it is not hard to see that in a $(2,1)$-category, having powers by $\mathcal{I}$ is equivalent to having powers by $\mathbf{2}$; natural transformations between functors $\mathbf{2} \to \K(X,Y)$ are all invertible. Therefore, if we show that categories internal to a $(2,1)$-category are groupoids, we can apply (Theorem 4.19, \cite{bourke2010codescent}) as desired. For $\X$ a category internal to $\K$, we have a $2$-cell $d_1 \Rightarrow d_0: X_1 \to X_0$ in $\K$; since $\K$ is a $(2,1)$-category, this is an isomorphism. Thus, by Yoneda, we obtain a morphism $(-)^{-1}: X_1 \to X_1$ which satisfies the requirement for $\X$ to be a groupoid internal to $\K$.
\end{proof}

We can therefore apply the work in the rest of this section abstractly to a $(2,1)$-category. 

\begin{thm}
    Let $\K$ be a $(2,1)$-category satisfying the axioms listed below. Then we have a $(2,1)$-equivalence $\K \simeq \GpdE$ with $\E$ a locally cartesian closed locos with coequalisers; consequently, there is a description of a type theoretic algebraic weak factorisation system on $\K$, and the right adjoint splitting of the comprehension category associated to this is a model of MLTT with strictly stable choices of $\Sigma$-, $\Pi$- and $\texttt{Id}$-types.   

    \begin{enumerate}
        \item $\mathcal{K}$ has pullbacks and powers by $\mathcal{I}$.
        \item $\mathcal{K}$ has codescent objects of groupoids internal to $\mathcal{K}$ whose source and target maps form a two-sided discrete fibration.
        \item Codescent morphisms are effective in $\mathcal{K}$.
        \item Discrete objects in $\mathcal{K}$ are projective, in the sense of Definition 4.13 of \cite{bourke2010codescent}.
        \item For every object $A \in \mathcal{K}$, there is a projective object $P \in \mathcal{K}$ and a codescent morphism $c: P \rightarrow A$.
        \item $\K$ is lextensive.
        \item Isofibrations in $\K$ are exponentiable.
        \item $\K$ has finite $2$-colimits. 
    \end{enumerate}
\end{thm}

\begin{proof}
   Let $\K$ be a $(2,1)$-category satisfying the above conditions. Conditions (1)-(5) imply that \cref{thm: bourke} can be applied and so $\K \simeq \GpdE$. Including (6)-(8) allows us to apply \cref{thm: final characterisation} and conclude that $\E$ is a locally cartesian closed locos with coequalisers. By \cref{thm: ttawfs}, we obtain a type theoretic algebraic weak factorisation system on $\K$ with the required properties. 
\end{proof}

\section{Examples and future directions}
\label{sec: examples}

We conclude by spelling out what \cref{thm: internal analogue} and \cref{thm: ttawfs} means for certain examples. To the author's knowledge, internal groupoid models of MLTT have not been considered previously, with the exception of some upcoming work by Awodey and Emmenegger, and Agwu. We conjecture that for $\E$ a locally cartesian closed locos with coequalisers, the effective model structure of \cite{gambino_henry_sattler_szumilo_2022} can be upgraded into a type theoretic algbraic weak factorisation system for simplicial objects in $\E$. This could give interesting variations of models of homotopy type theory in a similar way to the method used here. Indeed, it is possible to show that for such an $\E$, the effective model structure is cofibrantly generated and hence algebraic. We leave this work for further research. 

\subsection{Realisability models}

Fix $A$ a partial combinatory algebra (PCA) and consider the category of assemblies $\mathbf{Asm}_{A}$ over this PCA. This is a locally cartesian closed locos with coequalisers, and so $\Gpd(\mathbf{Asm}_A)$ gives a model of MLTT where types are realised by elements of the PCA. This gives a realisability model of MLTT, which is different from the approaches to $2$-dimensional effective considered by both Speight, and Awodey and Emmenegger (see \cite{awodey2025effective} for Awodey and Emmenegger's approach). 

Moreover, the full subcategory of \emph{modest} assemblies is also a locally cartesian closed locos with coequalisers, and so $\mathbf{Gpd}(\mathbf{Mod}_{A})$ gives us a modest realisability model of MLTT.

The inclusion of $\mathbf{Gpd}(\mathbf{Mod}_{A}) \hookrightarrow \mathbf{Gpd}(\mathbf{Asm}_{A})$ will be studied together with Sam Speight in future work; in fact, we can model type theory with a univalent universe of $0$-types using this observation. 

It is interesting to note that the effective topos \cite{hyland1982effective} is the exact completion of the category of assemblies over Kleene's first algebra. Moreover, there is a sense in which the map $\E \mapsto \GpdE$ is a kind of $(2,1)$-exact completion of a $1$-category \cite[Corollary 61]{bourke2014two}. Hence $\Gpd(\mathbf{Asm}_{K_1})$ could be described as some kind of effective $(2,1)$-topos. Moreover, this will not be a Grothendieck $(2,1)$-topos as its $1$-category of $0$-truncated objects (which is $\mathbf{Asm}_{K_1}$) is not a Grothendieck $1$-topos. 

We could also apply this approach to the effective topos itself. The effective topos $\mathbf{Eff}$ \cite{hyland1982effective} is a suitable environment for higher-order recursion theory, and is of interest to logicians interested in issues in computability. It is an elementary topos with natural numbers object. Therefore, we obtain an algebraic model structure on $\Cat(\mathbf{Eff})$. The associated internal model of MLTT gives a setting in which higher-order recursion is baked into the types in some way. We believe that this is related to the work of Anthony Agwu. 

\subsection{Arithmetic $\Pi$-pretoposes}\label{subsec: arithmetic}
Arithmetic $\Pi$-pretoposes were constructed to be the syntactic categories for extensional MLTT \cite{streicher1993investigations}. They are univalent universes of dependent type theory that satisfy axiom K and are closed under the empty type, unit type, sum types, dependent sum types,
propositional truncations, quotient sets, and parametrised natural numbers type. Therefore, the categorical constructions of our paper could equivalently be formulated in the language of extensional MLTT, which is the internal language on an arithmetic $\Pi$-pretopos--- see \cite[\S 3]{maietti2010joyal} for details. As a result, the statement and proof of \cref{thm: ttawfs} could be written using extensional MLTT. As a result, if we use extensional MLTT as our foundation for mathematics and construct the category of groupoids internal to this foundation, we still obtain a groupoid model of \emph{intensional} MLTT. This is a relative consistency result, akin to results in classical set theory which state that the consistency of ZF can be proven using ZF with a universe. Our result proves the consistency of intensional MLTT, given the consistency of extensional MLTT. 

One example of interest is taking $\E$ to be a model of Palmgren's constructive elementary theory of the category of sets \cite{Palmgren2012Constructivist}. such a thing is a constructive categorical model of Bishop's set theory. In such a setting, we get a constructive, intuitionistic version of Hofmann and Streicher's original groupoid model of MLTT. A concrete example of such a model of CETCS is given by the category of setoids \cite{emmenegger2020exact}.

\subsection{Elementary toposes with a natural numbers object}\label{subsec: elt topos}

Let $\E$ be an elementary topos with a natural numbers object. Any elementary topos with a natural numbers object is in particular a locally cartesian closed locos with coequalisers, and so both \cref{thm: algebraic model structure on CatE,thm: ttawfs} apply and so we obtain an algebraic model structure on $\CatE$ and also an algebraic model of MLTT given by the cloven isofibrations in $\GpdE.$ Below, we give multiple examples of categories that are elementary toposes with a natural numbers object.

\subsubsection{Constructive versions of classical results}\label{subsubsec: set}

An example of an elementary topos with a natural numbers object is given by $\E = \s.$ In this case, we obtain a constructive version of the classical model structure on $\Cat$; in fact, the axiom of choice is equivalent to the existence of the classical model structure on $\Cat$--- that is a model structure in which the fibrations are isofibrations and the cofibrations are injective on objects functors. If we do not assume the axiom of choice, then complemented inclusion on object functors are injective on object functors, but without assuming the law of the excluded middle injective on object functors are not complemented inclusion on object functors. Moreover, the weak equivalences of this model structure are functors which are fully faithful and essentially surjective on objects, which are actual equivalences of categories if and only if the axiom of choice holds \cite{freyd1990categories}. Therefore, assuming the axiom of choice, the (algebraic) model structure considered in this paper agrees with the classical model structure on $\Cat$, and so constructively this model structure is the correct one to consider. 

In this case, \cref{thm: ttawfs} provides an algebraic version of Hofmann and Streicher's groupoid model \cite{hofmann1998groupoid}; forgetting the algebraic structure, the display maps are given by isofibrations. Another algebraic version of Hofmann and Streicher's model is given by \cite[Theorem 4.5]{gambino2023models}. The algebraic weak factorisation system considered there has right algebras given by \emph{normal} isofibrations, which are cloven isofibration in which identities lift to identities. Assuming the axiom of choice, this agrees with our case; lifts can be chosen to be normal. However, constructively, there is a difference in these models of MLTT; in the type theory modelled by our result, identities are not decidable. 

In fact, isofibrations are equivalent to normal isofibrations in $\CatE$ if and only if injective morphisms in $\E$ are complemented inclusions  \cite[Corollary 3.12]{bourke20242}, and so identities are decidable in our model if and only if we assume the law of the excluded middle.

\subsubsection{For a presheaf category}

Let $\mathbb{C}$ be a small category. For $\E = [\C^{\op}, \s],$ we have $\CatE \cong [\C^{\op} , \Cat]$. To see this, note that $\X \in \Cat([\C^{\op}, \s])$ is 

\begin{equation*}
    \begin{tikzcd}
    ... \arrow[r] & X_1 \times_{X_0}X_1 \arrow[r, "m"] \arrow[r, "p_1", shift left = 4] \arrow[r, "p_2", shift right = 4, swap] & X_1 \arrow[r, shift left = 4, "d_1"] \arrow[r, shift right = 4, "d_0", swap] & X_0 \arrow[l, "i", swap]
\end{tikzcd}
\end{equation*}

with $X_0, X_1 \in [\C^{\op}, \s]$. Due to the coherence conditions for sources, targets, composition that means that for each $c \in \C$, we have a small category:

\begin{equation*}
    \begin{tikzcd}
    ... \arrow[r] & X_1(c) \times_{X_0(c)}X_1(c) \arrow[r, "m(c)"] \arrow[r, "p_1(c)", shift left = 4] \arrow[r, "p_2(c)", shift right = 4, swap] & X_1(c) \arrow[r, shift left = 4, "d_1(c)"] \arrow[r, shift right = 4, "d_0(c)", swap] & X_0(c) \arrow[l, "i", swap]
\end{tikzcd}
\end{equation*}

which assembles into a functor $\X: \C^{\op} \to \Cat.$

An internal isofibration $f: \X \to \Y$ is an internal isofibration if and only if \[\X(\mathcal{I}) \to \Y(\mathcal{I})\times_{Y_0} X_0\]

is a split epi in $[\C^{\op}, \s].$ A map in a presheaf category is split epi if and only if it is levelwise split epi, so $f$ is an internal isofibration if and only if for every $c\in \C^{\op}$

\[\X(\mathcal{I})(c) \to \Y(\mathcal{I})(c)\times_{Y_0(c)} X_0(c)\]

is a split epi in $\s$, where the equivalence $\Y(\mathcal{I} \times{Y_0} X_0)(c) \cong \Y(\mathcal{I})(c)\times_{Y_0(c)} X_0(c) $ is because limits are computed pointwise in presheaf categories. This happens if and only if $f(c): \X(c) \to \Y(c)$ is an isofibration in $\Cat.$ 

A similar calculation shows that $f: \X \to \Y$ is a weak equivalence if and only if $f(c): \X(c) \to \Y(c)$ is a weak equivalence in $\Cat.$ Therefore the structure of this model structure is calculated pointwise and is the projective model structure as described in \cite{hirschhorn2003model}. 

In this case our model structure coincides with the projective model structure lifted from the natural model structure on $\Cat$. In particular, for $\E = \sS,$ this agrees with Horel's levelwise model structure on $\Cat(\sS)$ given in Section 5.1 of \cite{horel2015model}. As such, we have given algebraic versions of these results. 

Any presheaf category has a natural numbers object given by the constant sheaf on the natural numbers. It is also an elementary topos, and therefore a locally cartesian closed locos with coequalisers, and so \cref{thm: ttawfs} can be applied and so $\Gpd([\C^{\op}, \s]) \cong [\C^{\op}, \Gpd]$ gives us an indexed model of MLTT. This could be of interest for further study.  

\subsubsection{For a Grothendieck topos}

A Grothendieck topos $$\begin{tikzcd}
    \E \arrow[r, hook, shift right=1.5] & \arrow[l, shift right =1.5] [\C^{op},\s] 
\end{tikzcd}$$
has a natural numbers object given by the sheafification of the natural numbers object in  $[\C^{op},\s]$.

In this case, our model structure does not usually agree with Joyal and Tierny's model structure on strong stacks \cite{joyal2006strong}; it agrees if and only if $\E$ satisfies the external axiom of choice; this occurs if and only if fully faithful and essentially surjective on objects functors are equivalences of categories. If the axiom of choice does not hold, then their weak equivalences are strictly contained in ours. In this case, our cofibrations are strictly contained in theirs as their cofibrations are injective on objects. Consequently, their fibrations are strictly contained in ours.

The associated model of MLTT given by $\GpdE$ can be thought of as having geometric types that are glued together from more simple indexed types.

For any Grothendieck topos satisfying the axiom of choice we have proven that strong stacks give a model of MLTT. Moreover, since Diaconescu's Theorem can be formulated in any topos, such a Grothendieck topos also satisfies the law of the excluded middle. Therefore, by \cite[Corollary 3.12]{bourke20242}, cloven isofibrations are equivalently describes as normal isofibrations and hence in the associated type theory, identities are decidable.

\subsection{$n$-fold Categories and internal $n$-fold Categories}
\label{subsec: nfold}

For $\E = \Cat,$ we recover the model structure on $\mathbf{DblCat}$ studied in (Section 8.2 \cite{fiore2008model}). As such, we obtain an algebraic version of this. 

We note that the work of this paper extends these results; if $\E$ is lextensive, locally finitely presentable and cartesian closed, then $\CatE$ is lextensive, locally finitely presentable and cartesian closed. Hence, by \cref{rem: locallyfp}, we obtain an algebraic model structure on $\mathbf{DblCat}(\E) : = \Cat(\CatE).$ This is a novel observation. Moreover, we see that $\mathbf{DblCat}(\E)$ is lextensive, locally finitely presentable and cartesian closed, so iterating this result gives us an algebraic model structure on $n\mathbf{FoldCat}(\E)$ for any $n \in \N$. In particular, we obtain (algebraic) model structures on $n\mathbf{FoldCat}$ for each $n$. This is different to the model structure on $n\mathbf{FoldCat}$ described in \cite{fiore2010thomason}, as it is different even when $n = 2,$ as noted in the introduction of \cite{fiore2010thomason}. Thus these examples are novel and of possible interest for further study.

However, it should be noted that by the association of our model structure with Lack's trivial model structure on a $2$-category, this in some sense takes a $2$-dimensional view of $\mathbf{nFoldCat}(\E),$ which perhaps neglects the higher-dimensional structure. 

We cannot use this recursion technique to get higher dimensional models of MLTT, since $\E$ locally cartesian closed does not imply $\CatE$ is locally cartesian closed. Indeed, for $\E = \s$, it is well-known that $\Cat$ is not locally cartesian closed.


\begin{thebibliography}{00}

\bibitem{anderson1978fibrations}
D.W. Anderson.
\newblock Fibrations and geometric realizations.
\newblock {\em Bulletin for the American Mathematical Society}, 84(5):765--788,
  1978.

\bibitem{awodey2023cartesian}
S. Awodey.
\newblock Cartesian cubical model categories.
\newblock {\href{https://arxiv.org/abs/2305.00893}{arXiv:2305.00893}}, 2023.

\bibitem{awodey2025effective}
S. Awodey, J. Emmenegger.
\newblock Toward the effective $2$-topos
\newblock {\href{https://arxiv.org/abs/2503.24279}{arXiv:2503.24279}}, 2025.

\bibitem{berg2018path} B. Van Den Berg.  \newblock Path categories and propositional identity types. \newblock \emph{ACM Transactions on
Computational Logic} (TOCL), \newblock 19(2):1–32, 2018.

\bibitem{berg2020univalent}
B. van den Berg. \newblock Univalent polymorphism. \newblock \emph{Annals of Pure and
Applied Logic}, \newblock 171(6):102793, 2020.


\bibitem{bezem2014model}
M.~Bezem, T.~Coquand, and S.~Huber.
\newblock A model of type theory in cubical sets.
\newblock In {\em 19th International Conference on Types for Proofs and
  Programs (TYPES 2013)}, volume~26, pages 107--128. Schloss Dagstuhl--Leibniz
  Zentrum fuer Informatik Wadern, Germany, 2014.

\bibitem{bourke2014two}
J. Bourke and R. Garner.
\newblock Two-dimensional regularity and exactness.
\newblock {\em Journal of Pure and Applied Algebra}, 218(7):1346--1371, 2014.

\bibitem{bourke20242}
J. Bourke and S. Lack.
\newblock On 2-categorical $\infty$-cosmoi.
\newblock {\em Journal of Pure and Applied Algebra}, 228(9):107661, 2024.

\bibitem{bourke2010codescent}
J. Bourke.
\newblock {\em Codescent objects in 2-dimensional universal algebra}.
\newblock PhD thesis, University of Sydney, 2010.

\bibitem{bunge1979stacks}
M. Bunge and R. Par{\'e}.
\newblock Stacks and equivalence of indexed categories.
\newblock {\em Cahiers de Topologie et G{\'e}om{\'e}trie Diff{\'e}rentielle
  Cat{\'e}goriques}, 20(4):373--399, 1979.

\bibitem{brown1973abstract}
K.S. Brown.
\newblock Abstract homotopy theory and generalized sheaf cohomology.
\newblock {\em Transactions of the American Mathematical Society},
  186:419--458, 1973.

\bibitem{barr2000toposes}
M. Barr and C. Wells.
\newblock {\em Toposes, triples, and theories}.
\newblock Springer-Verlag, 2000.

\bibitem{carboni1993introduction}
A. Carboni, S. Lack, and R.F.C Walters.
\newblock Introduction to extensive and distributive categories.
\newblock {\em Journal of Pure and Applied Algebra}, 84(2):145--158, 1993.

\bibitem{everaert2005model}
T. Everaert, R.W. Kieboom, and T. van~der Linden.
\newblock Model structures for homotopy of internal categories.
\newblock {\em Theory Appl. Categ}, 15(3):66--94, 2005.

\bibitem{emmenegger2020exact}
J. Emmenegger and E. Palmgren.
\newblock Exact completion and constructive theories of sets.
\newblock {\em The Journal of Symbolic Logic}, 85(2):563--584, 2020.

\bibitem{fiore2010thomason}
T.M. Fiore and S. Paoli.
\newblock A {Thomason} model structure on the category of small $n$-fold categories.
\newblock {\em Algebraic \& Geometric Topology}, 10(4):1933--2008, 2010.

\bibitem{fiore2008model}
T.M. Fiore, S. Paoli, and D. Pronk.
\newblock Model structures on the category of small double categories.
\newblock {\em Algebraic \& Geometric Topology}, 8(4):1855--1959, 2008.

\bibitem{freyd1990categories}
P.J. Freyd and A. Scedrov.
\newblock {\em Categories, Allegories}.
\newblock Elsevier, 1990.

\bibitem{garner2009understanding}
R. Garner.
\newblock Understanding the small object argument.
\newblock {\em Applied Categorical Structures}, 17(3):247--285, 2009.

\bibitem{gambino_henry_sattler_szumilo_2022}
N. Gambino, S. Henry, C. Sattler, and K. Szumi{\l}o.
\newblock The effective model structure and $\infty $ -groupoid objects.
\newblock {\em Forum of Mathematics, Sigma}, 10:e34, 2022.

\bibitem{gambino2023models}
N. Gambino and M.F. Larrea.
\newblock Models of {M}artin-{L}\"{o}f type theory from algebraic weak factorisation systems.
\newblock {\em The Journal of Symbolic Logic}, 88(1):242--289, 2023.

\bibitem{gambino2022constructive}
N. Gambino, C. Sattler, and K. Szumi{\l}o.
\newblock The constructive {K}an--{Q}uillen model structure: two new proofs.
\newblock {\em The Quarterly Journal of Mathematics}, 73(4):1307--1373, 2022.

\bibitem{grandis2006natural}
M. Grandis and W. Tholen.
\newblock Natural weak factorization systems.
\newblock {\em Archivum Mathematicum}, 42(4):397--408, 2006.

\bibitem{hirschhorn2003model}
P.S. Hirschhorn.
\newblock {\em Model categories and their localizations}.
\newblock Number~99. American Mathematical Soc., 2003.



\bibitem{Hughes2024Elementary}
C. Hughes and A. Miranda.
\newblock The elementary theory of the 2-category of small categories.
\newblock {\em Theory and Applications of Categories}, Vol. 43, 2025, No. 8, pp 196-242.

\bibitem{Hughes2024Colimits}
C. Hughes and A. Miranda.
\newblock Colimits of internal categories.
\newblock {\href{https://arxiv.org/abs/2501.17769}{arxiv:2501.17769}}, 2025.
\bibitem{horel2015model}
G. Horel.
\newblock A model structure on internal categories in simplicial sets.
\newblock {\em Theory and Applications of Categories}, 30(20):704--750, 2015.

\bibitem{hovey1998monoidal}
M. Hovey.
\newblock Monoidal model categories.
\newblock {\href{https://arxiv.org/abs/math/9803002}{	arXiv:math/9803002}}, 1998.

\bibitem{hofmann1998groupoid}
M. Hofmann and T. Streicher.
\newblock The groupoid interpretation of type theory.
\newblock {\em Twenty-five years of constructive type theory (Venice, 1995)},
  36:83--111, 1998.

\bibitem{hyland1982effective}
J.M.E Hyland.
\newblock The effective topos.
\newblock In {\em Studies in Logic and the Foundations of Mathematics}, volume
  110, pages 165--216. Elsevier, 1982.

\bibitem{johnstone2002sketches}
P.T. Johnstone.
\newblock {\em Sketches of an Elephant: A Topos Theory Compendium: Volume 2},
  volume~2.
\newblock Oxford University Press, 2002.

\bibitem{joyal2006strong}
A. Joyal and M. Tierney.
\newblock Strong stacks and classifying spaces.
\newblock In {\em Category Theory: Proceedings of the International Conference
  held in Como, Italy, July 22--28, 1990}, pages 213--236. Springer, 2006.

\bibitem{kelly1982basic}
M. Kelly.
\newblock {\em Basic concepts of enriched category theory}, volume Reprints in
  Theory and Applications of Categories.
\newblock Cambridge Univ. Press; Cambridge, MR0651714, 1982.

\bibitem{lack2006homotopy}
S. Lack.
\newblock Homotopy-theoretic aspects of 2-monads.
\newblock {\em Journal of Homotopy and Related Structures}, 2(2):229--260,
  2007.

\bibitem{Lamarche1991Proposal}
F. Lamarche.
\newblock A proposal about foundations i.
\newblock Manuscript. \url{https://www.cse.chalmers.se/~coquand/Proposal.pdf},
  1991.
\newblock [Online; accessed 21/03/2025].

\bibitem{maietti2010joyal}
M.E. Maietti.
\newblock Joyal's arithmetic universe as list-arithmetic pretopos.
\newblock {\em Theory \& Applications of Categories}, 24, 2010.

\bibitem{niefield2019internal}
S.B. Niefield and D.A. Pronk.
\newblock Internal groupoids and exponentiability.
\newblock {\em Cahiers de Topologie et G{\'e}om{\'e}trie Diff{\'e}rentielle
  Cat{\'e}goriques, LX}, 4, 2019.

\bibitem{Palmgren2012Constructivist} Erik Palmgren. \newblock Constructivist and structuralist foundations: Bishop’s and Lawvere’s theories of sets. \newblock \emph{Annals of Pure and Applied Logic,} \newblock 163(10):1384–1399, 2012.


\bibitem{quillen1967homotopical}
D.G. Quillen.
\newblock {\em Homotopical Algebra}.
\newblock Lecture notes in mathematics. Springer-Verlag, 1967.

\bibitem{rezk1996model}
C. Rezk.
\newblock A model category for categories.
\newblock {\em \href{https://ncatlab.org/nlab/files/Rezk_ModelCategoryForCategories.pdf}{Available from the author’s web page}}, 1996.

\bibitem{riehl2011algebraic}
E. Riehl.
\newblock {\em Algebraic model structures}.
\newblock The University of Chicago, 2011.

\bibitem{riehl2014categorical}
E.~Riehl.
\newblock {\em Categorical homotopy theory}, volume~24.
\newblock Cambridge University Press, 2014.

\bibitem{riehl2017type}
E. Riehl and M. Shulman.
\newblock A type theory for synthetic $\infty$-categories.
\newblock {\href{https://arxiv.org/abs/1705.07442}{arXiv:1705.07442}}, 2017.

\bibitem{streicher1993investigations}
T. Streicher.
\newblock Investigations into intensional type theory.
\newblock {\em Habilitiation Thesis, Ludwig Maximilian Universit{\"a}t},
  page~57, 1993.

\bibitem{program2013homotopy}
The {Univalent Foundations Program}.
\newblock {\em Homotopy Type Theory: Univalent Foundations of Mathematics}.
\newblock \url{https://homotopytypetheory.org/book}, Institute for Advanced
  Study, 2013.


\bibitem{voevodsky2006very}
V. Voevodsky.
\newblock A very short note on the homotopy $\lambda$-calculus, 2006.

\end{thebibliography}
\end{document}